\documentclass[10pt, peerreview, draftclsnofoot]{IEEEtran}
\pdfoutput=1
\usepackage{cite, graphicx, subfigure,url,setspace,algorithm}
\usepackage{amsmath,amssymb}
\usepackage{bm}
\usepackage{subfigure}
\usepackage{CJK}

\newtheorem{proposition} {Proposition}
\newtheorem{remark}{Remark}

\usepackage{epstopdf}
\begin{document}

\title{ Multi-Target Tracking Using A Randomized Hypothesis Generation Technique}

\author{
\IEEEauthorblockN{W. Faber and S. Chakravorty}
\IEEEauthorblockA{Department of Aerospace Engineering \\ Texas
  A\&M University \\ College Station, TX} 
  \and
\IEEEauthorblockN{Islam I. Hussein}\IEEEauthorblockA{Applied Defense Solutions \\ Columbia, MD}}
\maketitle

\setlength{\parskip}{0ex}

\begin{abstract}
 In this paper, we present a randomized version of the finite set statistics (FISST) Bayesian recursions for multi-object tracking problems. We propose a hypothesis level derivation of the FISST equations that shows that the multi-object tracking problem may be considered as a finite state space Bayesian filtering problem, albeit with a growing state space. We further show that the FISST and Multi-Hypothesis Tracking (MHT) methods for multi-target tracking are essentially the same. We propose a randomized scheme, termed randomized FISST (R-FISST), where we sample the highly likely hypotheses using Markov Chain Monte Carlo (MCMC) methods which allows us to keep the problem computationally tractable. We apply the R-FISST technique to a fifty-object birth and death Space Situational Awareness (SSA) tracking and detection problem. We also compare the R-FISST technique to the Hypothesis Oriented Multiple Hypothesis Tracking (HOMHT) method using an SSA example. 
\end{abstract}

\section{Introduction}
In this paper, we present a randomized approach to approximate the full Bayesian recursions involved in solving the Finite Set Statistics (FISST) based approach to the problem of multi-object tracking and detection.  We show that the FISST recursions can essentially be considered as a discrete state space Bayesian filtering problem on ``Hypothesis Space" with the only input from the continuous problem coming in terms of the likelihood values of the different hypotheses. The number of objects is implicit in this technique and can be a random variable. The "Hypothesis Space" perspective allows us to develop a randomized version of the FISST recursions where we sample the highly likely children hypotheses using a Markov Chain Monte Carlo (MCMC) technique thereby allowing us to keep the total number of possible hypotheses under control, and thus, allows for a computationally tractable implementation of the FISST equations, which otherwise grows at an exponential rate, and thus, can quickly lead to the problem becoming intractable. The method is applied to a fifty-object SSA tracking and detection problem that has an unfixed number of objects throughout the simulation. The method is then compared to that of a well-known tracking method, HOMHT, using a fifteen-object SSA example and a variety of data association scenarios. 

In the last 20 years, the theory of FISST-based multi-object detection and tracking has been developed based on the mathematical theory of finite set statistics \cite{Mahler:97,Mahler:07}.  The greatest challenge in implementing FISST in real-time, which is critical to any viable SSA solution, is computational burden. The first-moment approximation of FISST  is known as the Probability Hypothesis Density (PHD) approach \cite{Mahler:07,Vo06}. The PHD has been proposed as a computationally tractable approach to applying FISST. The PHD filter essentially finds the density of the probability of an object being at a given location, and thus, can provide information about the number of objects (integral of the PHD over the region of interest) and likely location of the objects (the peaks of the PHD). The PHD can further employ a Gaussian Mixture (GM) or a particle filter approximation to reduce the computational burden (by removing the need to discretize the state space). This comes at the expense of approximating the general FISST pdf with its first-moments \cite{Vo06, Vo_SMC, Vo_CPHD, Clark_PHD}. The PHD filter does not attempt to solve the full FISST recursions, in particular, by considering the PHD, the filter gets rid of the data association problem inherent in these problems. In other previous work, a GM approximation was applied, not to the first-moment of the FISST pdfs, but to the original full propagation and update equations derived from FISST \cite{Hussein:12GNC,Hussein:12b}. This eliminates any information loss associated with using the first-moment PHD approximation, while at the same time increasing the computational tractability of the multi-object FISST pdfs. This approach is similar in spirit to the concept of the ``para-Gaussian" pdf that was described in \cite{Mahler:99}.

In this paper, in contrast, we introduce a hypothesis level derivation of the FISST equations that makes it clear as to how the full FISST recursions can be implemented without any approximation other than the approximations inherent in the underlying tracking filter, such as an extended  Kalman filter. We introduce a simplified model for the birth and death process that allows for only one birth or one death in any time interval, thereby controlling the number of birth and death hypotheses while still being able to converge to the correct hypothesis regarding the total number of objects given the birth objects remain in the field of view for long enough. Further, in order to ensure the computational tractability of the resulting equations, we introduce an MCMC based hypothesis selection scheme resulting in the Randomized FISST (R-FISST) approach that is able to scale the FISST recursions to large scale problems. We call our method R-FISST, since as in FISST, the hypotheses in our method have varying number of objects, and in essence, give a probabilistic description of the random finite set representing the multi-target probability distribution. We also formally show the equivalence of the two methods in the appendix.

There are also non-FISST based approaches to multi-hypothesis tracking (MHT) such as the Hypothesis Oriented MHT (HOMHT) \cite{HOMHT, BarShalom1, BarShalom2, mcmcda}, and the track oriented MHT (TOMHT) techniques \cite{TOMHT}.  The MHT techniques can be divided into single-scan and multi-scan methods depending on whether the method uses data from previous times to distinguish the tracks \cite{BarShalom1, jpda, mcmcda}. The single-scan (recursive) methods such as joint probabilistic data association (JPDA) \cite{jpda, mcmcda} typically make the assumption that the tracks are independent which is not necessarily true. The multi-scan methods such as TOMHT \cite{TOMHT, mcmcda} are not recursive. The primary challenge in these methods is  the management of the various different hypotheses related to the tracks which TOMHT does using an efficient tree structure, and the MCMCDA, and other related tracking techniques \cite{mcmcda1, mcmcda2, mcmcda3}, do through the use of MCMC methods in sampling the data associations. We also use MCMC to sample children hypotheses given the parent hypothesis, however, our approach is a truly recursive technique which does not assume track independence as the above mentioned single scan methods. We essentially do an efficient management of the growing number of hypotheses at every generation through the judicious use of MCMC.
A primary contribution of this paper is to show that the MHT technique and the FISST technique are essentially the same, modulo the set-theoretic representation of multi-target pdfs in FISST (which, however, does not provide any extra information). This is made possible through the hypothesis level derivation of the tracking equations that considers the full hybrid state of the problem. This allows us to identify the critical structure inherent in the FISST recursions that enables us to unify the two approaches.
%Also, comparisons and research on implementing the randomized technique with a HOMHT framework as well as a more detailed discussion of the randomized technique' and what distinguishes it from other randomized methods, such as MCMCDA, shall be presented in the proceedings of the 2015 AAS/AIAA Astrodynamics Specialist Conference to be held at Vail, CO, August 9-13 of 2015.
The contributions of this paper are as follows:
\begin{enumerate}
\item{ We introduce a hypothesis level derivation of the FISST recursions and show its equivalence to the standard FISST recursions. This allows us to implement the full FISST recursions as opposed to first moment approximations of the recursions such as PHD.}
\item{We unify the FISST and MHT methodologies based on the hypothesis level derivation.}
\item{We propose a randomized MCMC based hypothesis generation technique, called RFISST, that allows us to keep the FISST recursions computationally tractable.}
\item{We apply our methods to the SSA problem and perform a detailed comparison of our technique with the HOMHT technique.}
\end{enumerate}

The rest of the paper is organized as follows. In Section II, we introduce the hypothesis level derivation of the FISST equations. In Section III, we introduce the MCMC based randomized hypothesis selection technique that results in the RFISST algorithm. In Section IV, we show an application of the RFISST technique to a fifty-object birth death SSA scenario. Section V discusses the comparison with HOMHT. In the appendix, we show the equivalence of the FISST and the hypothesis level FISST equations. This paper is an expanded and unified version of the references \cite{Faber1},\cite{Faber2}.   %A related paper,\cite{Faber1}, was presented at the International Conference of Information Fusion. This paper extends that paper with a new MCMC data association scheme and a full comparison of the methodology with the HOMHT technique. For the sake of the paper being self-contained, the next two sections detailing the hypothesis level derivation of the FISST equations are reproduced from reference \cite{Faber2}

\section{A Hypothesis based Derivation of the FISST equations}
In this section, we shall  frame the multi-object tracking equations at the discrete hypothesis level ( which we believe are equivalent to the FISST equations) which then shows clearly as to how the full FISST recursions may be implemented. The derivation below assumes that the number of measurements is always less than the number of objects, which is typically the case in the SSA problem.  We never explicitly  account for the number of objects, since given a hypothesis, the number of objects and their probability density functions (pdf) are fixed, which  allows us to derive the results without having  to consider the random finite set (RFS) theory underlying FISST.  Albeit the equations derived are not as general as the FISST equations, in particular, the birth and death models employed here are quite simple, we believe that the level of generality is sufficient for the SSA problem that is our application. 

\subsection{Framing FISST at the Hypothesis Level}
We consider first the case when the number of objects is fixed, which we shall then generalize to the case when the number of objects is variable, i.e, there is birth and death in the object population. Assume that the number of objects is $M$, and each object state resides in $\Re^N$. Consider some time instant $t-1$, and the data available for the multi-object tracking problem till the current time $\mathcal{F}^{t-1}$. Let $H_i$ denote the $i^{th}$ hypothesis at time $t-1$, and let $\{X\}$ denote the underlying continuous state. For instance, given the $N-$ object hypothesis, the underlying state space would be $\{X\} = \{X_1, X_2,\cdots X_N\}$ where $X_j$ denotes the state of the $j^{th}$ object under hypothesis $H_i$ and resides in $\Re^N$. Let $p(\{X\}, i/ \mathcal{F}^{t-1})$ denote the joint distribution of the state-hypothesis pair after time $t-1$. Using the rule of conditional probability:
\begin{align}
p(\{X\}, i/ \mathcal{F}^{t-1}) = \underbrace{p(\{X\}/ i, \mathcal{F}^{t-1})}_{\mbox{MT-pdf underlying} H_i} \underbrace{p(i/ \mathcal{F}^{t-1})
}_{w_i = \mbox{prob. of} H_i},
\end{align} 
where MT-pdf is the multi-object pdf underlying a hypothesis. Given the hypothesis, the MT-pdf is a product of independent individual pdfs underlying the objects, i.e.,
\begin{align}
p(\{X\}/ i, \mathcal{F}^{t-1}) = \prod_{k=1}^M p_k(x_k),
\end{align}
where $p_k(.)$ is the pdf of the $k^{th}$ object.
\begin{remark}
In random finite set theory, the arguments of the MT-pdf $\{x_1, x_2 \cdots x_M\}$ above are interchangeable and thus, the MT-pdf is represented as:
\begin{align}
p(\{X\}/ i, \mathcal{F}^{t-1}) = \sum_{\bar{\sigma}}\prod_{k=1}^M p_{\sigma_k}(x_k),
\end{align}
where $\bar{\sigma} = \{\sigma_1, \sigma_2 \dots \sigma_M\}$ represents all possible permutations of the indices $\{1,2 \cdots M\}$. Hence, in any integration involving such a set of indices, a normalizing factor of $\frac{1}{M!}$ is used. In our case, we explicitly assign the index $x_k$ to the target $k$ or more precisely, the $k^{th}$ component of the MT-pdf, $p_k(.)$. Note that such an assignment is always possible and there is no information loss in such a representation. Moreover, at the expense of a bit more bookkeeping, this allows us to keep track of the labels of the different components of our multi-target pdfs. Please also see the appendix where we show the equivalence of the hypothesis level equations derived here and the FISST recursions.
\end{remark} 

Next, we consider the prediction step between measurements. Each hypothesis $H_i$ splits into $A_M$ children hypotheses, and let us denote the $j^{th}$ child hypothesis as $H_{ij}$. The children hypotheses correspond to the different data associations possible given a measurement of size $m$, i.e., $m$ returns, and
\begin{align}\label{n_da}
A_M = \sum_{n=0}^{min(m,M)} {M \choose n} {m \choose n} n!.
\end{align}
We want to note here that this is a pseudo-prediction step since we assume that we know the size of the return $m$. However, it allows us to fit the MT-tracking method nicely into a typical filtering framework.
Using the rules of total and conditional probability, it follows that the predicted multi-object pdf in terms of the children hypotheses is:
\begin{align} \label{FISST_pred}
p^-(\{X\}, (i,j)/ \mathcal{F}^{t-1}) = %\nonumber\\
\int  p(\{X\}, (i,j)/ \{X'\}, i)p(\{X'\}, i/ \mathcal{F}^{t-1})d\{X'\} = \nonumber\\
\underbrace{\int p(\{X\}/(i,j), \{X'\})p(\{X'\}/ i, \mathcal{F}^{t-1}) d\{X'\}}_{p^-(\{X\}/ (i,j), \mathcal{F}^{t-1})}%\nonumber\\
\underbrace{p(j/ i)}_{p_{ij}} 
\underbrace{p(i/ \mathcal{F}^{t-1})}_{w_i},
\end{align}
where $p^{-}(., (i,j)/ \mathcal{F}^{t-1})$ is the joint distribution of the state and hypothesis pairs before the measurement at time $t$. We have used the fact that $p((i,j)/ \{X'\}, i) = p(j/i) = p_{ij}$, and $p_{ij}$ is the transition probability of going from the parent $i$ to the child $j$ and $w_i$ is the probability of the parent hypothesis $H_i$.  Let $p_k(x_k/ x_k')$ denote the transition density function of the $k^{th}$ object. Expanding the predicted MT-pdf, we obtain:
\begin{align}
p^-&(\{X\}/ (i,j), \mathcal{F}^{t-1}) \equiv %\nonumber\\
\int p(\{X\}/ (i,j), \{X'\})p(\{X'\}/ (i), \mathcal{F}^{t-1})d\{X'\}, \end{align}
where
\begin{align}
\label{pred_nbd}
p(\{X\}  / (i,j), \{X'\}) \equiv \prod_{k=1}^M p_k(x_k/ x_k') %\nonumber\\
\int  p(\{X\}/ (i,j), \{X'\})p(\{X'\}/ (i), \mathcal{F}^{t-1})d\{X'\}  \nonumber\\
\equiv\int \prod_k p_k(x_k/ x_k') \prod_{k'} p_{k'}(x_k') dx_1'\cdots dx_M'% \nonumber\\
= \prod_ k \int p_k(x_k/ x_k') p_k(x_k') dx_k' = \prod_k p_k^-(x_k),
\end{align}
where $p_k^-(x_k)$ is the prediction of the $k^{th}$ object pdf underlying the hypothesis $H_{ij}$. 
\begin{remark}
Eq. \ref{FISST_pred} has a particularly nice hybrid structure: note that the first factor is the multi-object continuous pdf underlying the child hypothesis $H_{ij}$, while the second factor $p_{ij} w_i$ is the predicted weight of the hypothesis $H_{ij}$.
For the no birth and death case, in the absence of any a priori information regarding the sensor,  all $p_{ij}$ are equal to $\frac{1}{A_M}$, where recall that $A_M$ is the total number of data associations possible (Eq. \ref{n_da}). However, if a priori information, for instance, in terms of a probability of detection $p_D$ is available, then:
\begin{align} \label{hyp_prior_wt}
p_{ij} = \frac{p_D^k (1-p_D)^{M-k}}{{m\choose k} k!},
\end{align}
if $ij$ is a data association in which $k$ of the $M$ targets are associated to measurements. The ${m\choose k} k!$ factor is required so that $p_{ij}$ is a valid probability distribution , i.e., $\sum_j p_{ij} = 1$. To see this, note that there are exactly ${m\choose k} k!$ ways that $m$ measurements can be assigned to $k$ targets, and there are $M\choose k$ ways of choosing the $k$ targets to associate to the $k$ chosen measurements. We are assuming here that $M< m$ which is almost always the case. Thus:
\begin{align}
\sum_j p_{ij} = \sum_{k=0}^M {M \choose k} \frac{p_D^k(1-p_D)^{M-k}}{{m\choose k}k!} {m\choose k}k! = 1.
\end{align}
Note that the MT-pdf underlying $H_{ij}$ is simply the product of the predicted individual object pdf, and in the case of no birth and death, it is the same for all children hypothesis $H_{ij}$. 
\end{remark}

Given the prediction step above, let us consider the update step given the measurements $\{Z_t\} = \{z_{1,t},\cdots z_{m,t}\}$, where there are $m$ measurement returns. We would like to update the weights of each of the multi-object hypotheses to obtain $p(\{X\}, (i,j)/ \{Z_t\},\mathcal{F}^{t-1})$ by incorporating the measurement $\{Z_t\}$. Using Bayes rule:
%\begin{widetext}
\begin{align}
p(&\{X\}, (i,j)/ \{Z_t\},\mathcal{F}^{t-1}) = %\nonumber\\
\eta p(\{Z_t\}/ \{X\}, (i,j)) p^-(\{X\}, (i,j)/ \mathcal{F}^{t-1}), \nonumber\end{align}
where
\begin{align}
\eta = 
 \sum_{i',j'} \int p(\{Z_t\}/ \{X'\}, (i', j')) p^-(\{X'\}, (i',j')/ \mathcal{F}^{t-1}) d\{X'\},\nonumber
\end{align}
where the MT-likelihood function $p(\{Z_t\}/ \{X\}, (i,j))$ and the Bayes normalizing factor $ \int p(\{Z_t\}/ \{X'\}, (i', j')) p^-(\{X'\}, (i',j')/ \mathcal{F}^{t-1}) d\{X'\}$ are defined in Eqs. \ref{MT-likelihood} and \ref{MT-Bayes_factor} below.
%\end{widetext}
Using the prediction equation \ref{FISST_pred}, it follows that:
\begin{align}
\underbrace{p(\{X\}, (i,j)/ \{Z_t\}, \mathcal{F}^{t-1})}_{p(\{X\}, (i,j)/ \mathcal{F}^t)} = %\nonumber\\
\hspace{2ex}\eta p(\{Z_t\}/ {X}, (i,j)) p^-(\{X\}/ (i,j), \mathcal{F}^{t-1})p_{ij}w_i.
\end{align}
We may then factor the above equation as follows:
\begin{align}\label{FISST_M}
p(\{X\}, (i,j)/ \mathcal{F}^t) = %\nonumber\\
\frac{p(\{Z_t\}/ \{X\}, (i,j))p^-(\{X\}/ (i,j), \mathcal{F}^{t-1})}{l_{ij}}
\frac{l_{ij}\overbrace{p_{ij}w_i}^{w_{ij}}}{\sum_{i',j'} l_{i',j'} \underbrace{p_{i'j'}w_{i'}}_{w_{i'j'}}},
\end{align}
where
\begin{align}
l_{ij} = \int p(\{Z_t\}/ \{X'\}, (i,j))p^-(\{X'\}/ (i,j), \mathcal{F}^{t-1}) d\{X'\}.
\end{align}
Note that $l_{ij}$ is likelihood of the data $\{Z_t\}$ given the multi-object pdf underlying hypothesis $H_{ij}$, and the particular data association that is encoded in the hypothesis. 
\begin{remark}
It behooves us to understand the updated pdf underlying the child hypothesis $H_{ij}$, the first factor on the right hand side of Eq. \ref{FISST_M}. Let $p_D$ denote the probability of detection of a object given that it is in the field of view (FOV) of the monitoring sensor(s). Let $g(z)$ denote the probability that the observation $z$ arises from a clutter source. Let $H_i$ denote an $M-$object hypothesis with object states $\{X\}=\{X_1,\cdots X_M\}$ governed by the pdfs $p_1(x_1), \cdots p_M(x_M)$. Let the child hypothesis $H_{ij}$ correspond to the following data association hypothesis: $z_1 \rightarrow X_{j_1}, \cdots z_m \rightarrow X_{j_m}$. Then, we define the MT-likelihood function:
\begin{align} \label{MT-likelihood}
p(\{Z_t\}/ \{X\}, (i,j)) \equiv %\nonumber\\
p(\{z_1\cdots z_m\}/ \{X_1 = x_1, \cdots X_M = x_M\}, (i,j)) %\nonumber\\
=[ \prod_{k=1}^m p(z_k/ X_{j_k} = x_{j_k})], 
\end{align}
where $p(z_k/ X_{j_k} = x_{j_k})$ is simply the single object observation likelihood function for the sensor.
Thus, 
\begin{align}
p(\{Z_t\}/\{X\}, (i,j))p^-(\{X\}/(i,j), \mathcal{F}^{t-1}) =%\nonumber\\
[\prod_{k=1}^m p(z_k/ X_{j_k}= x_{j_k})p_{j_k}^-(x_{j_k}) ] [\prod_{l\neq j_k}p_l^-(x_l)],
\end{align}
where $l \neq j_k$ denotes all objects $X_l$ that are not associated with a measurement under hypothesis $H_{ij}$.
Further, defining the MT-Bayes factor as:
\begin{align}\label{MT-Bayes_factor}
l_{ij}= \int p(\{Z_t\}/ \{X'\}, (i,j))p^-(\{X'\}/ (i,j), \mathcal{F}^{t-1}) d\{X'\}  \nonumber\\
\equiv\int [\prod_{k=1}^m p(z_k/ X_{j_k}= x'_{j_k}) p_{j_k}^-(x'_{j_k})] \times%\nonumber\\
[\prod_{l\neq j_k}p_l^-(x'_l)]dx_1'..dx_M' \nonumber\\
= [\prod_{k=1}^m (\int p(z_k/ X_{j_j}= x'_{j_k})p_{j_k}^-(x_{j_k}') dx'_{j_k})]\times%\nonumber\\
\hspace{2ex}[\prod_{l\neq j_k}\int p_l^-(x_l')dx_l']\nonumber\\
=  \prod_{k=1}^m p(z_k/ X_{j_k}),
\end{align}
where $p(z_k/ X_{j_k})\equiv \int p(z_k/ X_{j_k}= x'_{j_k}) p_{j_k}^-(x'_{j_k}) dx'_{j_k}$.
Hence,
\begin{align}
\frac{p(\{Z_t\}/\{X\}, (i,j))p^-(\{X\}/(i,j), \mathcal{F}^{t-1})}{\int p(\{Z_t\}/ \{X'\}, (i,j))p^-(\{X'\}/ (i,j), \mathcal{F}^{t-1}) d\{X'\} }%\nonumber\\
= \frac{\prod_{k=1}^m p(z_k/ X_{j_k} = x_{j_k})p^-_{j_k}(x_{j_k})}{\prod_{k=1}^m \int p(z_k/ X_{j_k} = x'_{j_k})p^-_{j_k}(x'_{j_k})dx'_{j_k}} \times%\nonumber\\
\prod_{l\neq j_k} p_l^-(x_l) \nonumber\\
= \prod_{k=1}^m p_{j_k}(x_{j_k}/ z_k) \times \prod_{l\neq j_k}p_l^-(x_l),
\end{align}
where $p_{j_k}(x_{j_k}/ z_k)$ denotes the updated object pdf of $X_{j_k}$ using the observation $z_k$ and the predicted prior pdf $p_{j_k}^-(x_{j_k})$, and $p_l^-(x_l)$ is the predicted prior pdf of $X_l$ whenever $l \neq j_k$, i.e., the pdf of object $X_l$ is not updated with any measurement. In the above, we have assumed that all the measurements are assigned to objects, however, some of the measurements can also be assigned to clutter, in which case, the object pdfs are updated exactly as above, i.e., all objects' predicted prior pdfs associated with data are updated while the unassociated objects' predicted priors are not updated,  except now the likelihoods $l_{ij}$of the children hypothesis $H_{ij}$ are given by:
\begin{align}\label{hyp_likelihood}
l_{ij} = \prod_{i=1}^m p(z_i/ X_{j_i}),
\end{align}
where
\begin{align}
p(z_i/ X_{j_i}) = 
\begin{cases}
    \int p(z_i/ x) p_{j_i}(x) dx & \text{if } X_{j_i} \in \mathcal{T}\\
    g(z_i)  & \text{if }  X_{j_i} \in \mathcal{C} 
    \end{cases}
\end{align}
where $\mathcal{T}$ is the set of all objects and $\mathcal{C}$ is clutter, $m'$ is the number of objects associated to measurements, and the above equation implies that the measurement $z_i$ was associated to clutter if $X_{j_i} \in \mathcal{C}$.
\end{remark}
\begin{remark}
The recursive equation \ref{FISST_M} above has a particularly  nice factored hybrid form.  The first factor is just a continuous multi-object pdf that is obtained by updating the predicted multi-object pdf obtained by associating the measurements in $\{Z_t\}$ to objects according to the data association underlying $H_{ij}$. The second factor corresponds to the update of the discrete hypothesis weights.
\end{remark}
\begin{remark}
Given that there is an efficient way to predict/ update the multi-object pdfs underlying the different hypotheses, Eq. \ref{FISST_M} actually shows that the FISST recursions may essentially be treated as a purely discrete problem living in the ``Hypothesis level" space. The "hypothesis level" weights are updated based on the likelihoods $l_{ij}$ which is determined by the continuous pdf underlying $H_{ij}$. Also,  the continuous pdf prediction and updates are independent of the hypothesis level prediction and updates, i.e, the hypothesis probabilities do no affect the multi-object pdfs underlying the hypotheses.
\end{remark} 
Thus, given that the likelihoods of different hypothesis $l_{ij}$ arise from the underlying multi-object pdf and the encoded data association in the hypotheses $H_{ij}$, the FISST updates can be written purely at the hypothesis level as follows:
\begin{align}\label{FISST_M'}
w_{ij} := \frac{l_{ij}w_{ij}}{\sum_{i',j'} l_{i'j'}w_{i'j'}},
\end{align}
where $w_{ij} = p_{ij}w_i$. Thus, we can see that the FISST update has a particularly simple Bayesian recursive form when viewed at the discrete hypothesis level, given that the multi-object pdfs underlying the hypotheses $H_{ij}$ are tracked using some suitable method.  We can summarize the above development of the Bayesian recursion for multi-object tracking as follows:
\begin{proposition}
Given an $M-$object hypothesis $H_i$, and its children hypotheses $H_{ij}$, that correspond to the data associations $\{z_i \rightarrow X_{j_i}\}$, the joint MT-density, hypothesis weight update equation is:
\begin{align}
p(\{X\}, (i,j)/ \mathcal{F}^t) = p(\{X\}/(i,j), \mathcal{F}^t)\frac{w_{ij}l_{ij}}{\sum_{i',j'} w_{i'j'}l_{i'j'}}, \nonumber
\end{align}
where $w_{ij} = p_{ij}w_i$, $l_{ij}$ is given by Eq. \ref{hyp_likelihood}, and the MT-pdf underlying $H_{ij}$:
\begin{align}
p(\{X\}/(i,j), \mathcal{F}^t) = \prod_{k=1}^m p_{j_k}(x_{j_k}/ z_k) \prod_{l\neq j_k} p_l^-(x_l), \nonumber
\end{align}
where $p_{j_k}(X_{j_k}/ z_{j_k})$ denotes the predicted prior of object $X_{j_k}$, $p_{j_k}^-(x_k)$, updated by the observation $z_{j_k}$, and $p_l^-(x_l)$ is the predicted prior for all objects $X_l$ that are not associated.
\end{proposition} 

We may renumber our hypothesis $H_{ij}$ into a parent of the next generation of hypothesis through a suitable map $F((i,j))$ that maps every pair $(i,j)$ into a unique positive integer $i'$, and start the recursive procedure again. However, the trouble is that the number of hypotheses grows combinatorially at every time step since at every step the number of hypotheses grow by the factor $A_M$ (Eq. \ref{n_da}), and thus, the above recursions can quickly get intractable.

\subsection{Relationship to MHT} 
The equations derived above are essentially the same equations as those derived in the MHT framework except for the ${m \choose k} k!$ factor in each hypothesis weight that is needed to normalize the MT-likelihood function similar to the case of a standard likelihood function. More specifically, consider Eqs. \ref{hyp_prior_wt} and \ref{MT-likelihood}. These two equations assure us that:
\begin{align}
\int p(\{z_1,.., z_m\}/\{x_1,..x_M\}, i) dz_1..dz_m %\nonumber\\
= \int \sum_{j} p_{ij} \prod_{l=1}^m p(z_l/ X_{j_l} = x_{j_l}) dz_1..dz_m 
 = 1,
\end{align}
for any MT-state $\{x_1,..x_M\}$, given that the parent hypothesis $i$ is an M-target hypothesis,
i.e., the total likelihood of all possible children hypothesis resulting from a parent hypothesis normalizes to unity, something that is typically required of the likelihood function $p(z/x)$ in a standard tracking/ filtering problem. The MHT likelihood for a child hypothesis $j$ of parent $i$ is of the form:
\begin{align}
\eta_{ij}^{MHT}(z_1,..z_m) = p_D^k (1-p_D)^{M-k} \prod_{l=1}^m p(z_l/ \hat{x}_{j_l}),
\end{align}
where $\hat{x}_{j_l}$ is the mean of the pdf of the target $X_{j_l}$ or $p(z_k/ \hat{x}_{j_l}) = g(z_l)$ if measurement  $z_l$ is associated to clutter. The equivalent hypothesis likelihood in the hypothesis level FISST (H-FISST) derivation is:
\begin{align}
\eta_{ij}^{HFISST}(z_1,..z_m) = p_{ij} l_{ij}=  \frac{{p_D}^k (1-p_D)^{M-k}}{{m \choose k}k!} \prod_{l=1}^m p(z_l/ X_{j_l}),
\end{align}
where the terms under the product in the above equation have been defined in Eq. \ref{hyp_likelihood}. Thus, it may be seen that the main difference in the two likelihoods is the factor ${m \choose k} k!$ and the fact that $p(z_l/ \hat{x}_{j_l})$ is an approximation of $p(z_l/ X_{j_l})$ for observations that have been associated to a target. Thus, in general, the weights of the different hypotheses will be different due to the MT-likelihood normalization required in H-FISST and the approximation of the true likelihood of target-observation associations in MHT, however, the MT-pdfs underlying the different hypotheses in the two methods are exactly the same. The normalization of the likelihood is necessary from a probabilistic perspective since otherwise the distribution on the filtered state pdf, resulting from all possible observations, does not constitute a probability distribution, i.e., it does not add up to unity. This can easily be seen for the case of a standard filtering problem which carries over to the mutli-target tracking problem. Let the filtered state pdf, the belief state be denoted by $b(x)$. Suppose that the likelihood function $\int p(z/x) dz \neq 1$. Consider the distribution on the future belief state $b'(x)$. This is given by:
\begin{align}
p(b'/b) = \int_z p(b'/z,b) p(z/b) dz, \mbox{where} %\nonumber\\
p(z/b) = \int p(z/x) b(x) dx.
\end{align}
Note that if $\int p(z/x)dz \neq 1$ then $\int p(z/b) dz \neq 1$. Hence,
\begin{align}
\int p(b'/b) db' = \int\int p(b'/z,b) p(z/b) dz db' %\nonumber\\
=  \int p(z/b) dz \neq 1.
\end{align}
We know that the filtered pdf (the belief process) has to evolve according to a Markov chain \cite{bertsekas1, kumar1} but the above development shows that the evolution equation violates the requirement that the transition probability of a Markov chain needs to be a probability distribution over all future states,  if the likelihood does not normalize to unity.
 \\

Furthermore, the proposed hybrid derivation (in that it includes both the continuous and discrete parts of the problem) as opposed to MHT which is a purely discrete derivation at the hypothesis level \cite{HOMHT},  reveals the critical hybrid structure (Eq. \ref{FISST_M}) inherent to  multi-target tracking problems, and which, in turn  allows us to unify the HFISST development with the FISST based formulation of the multi-target tracking problem, and thus, allows for the unification of FISST and MHT: methods that have thus far been thought to be different from each other (please see the appendix for more details) modulo the difference in the hypothesis weights due to the MT-likelihood normalization required in the FISST formulation. 

\subsection{Incorporating Birth and Death in Hypothesis level FISST}
The development thus far in this section has assumed (implicitly) that there are a fixed and known number of objects. However, this is not necessarily true since new objects can arrive while old objects can die. Thus, we have to incorporate the possibility of the birth and death of objects. In the following, we show that this can be done in quite a straightforward fashion using Eqs. \ref{FISST_pred},  \ref{FISST_M} and \ref{FISST_M'}. 

Let $\alpha$ denote the birth probability of a new object being spawned and $\beta$ denote the probability that an object dies in between two measurements. We will assume that $\alpha^2, \beta^2 \approx 0$. This assumption implies that exactly one birth or one death is possible in between measurement updates. 
Consider the time instant $t$, and consider an $M$-object hypothesis at time $t$, $H_i$.  
Depending on the time $t$, let us assume that there can be $M_t^b$ birth hypotheses and $M_t^d$ death hypothesis corresponding to one of $M_t^b$ objects being spawned or one of $M_t^d$ objects dying. In particular, for the SSA problem, we can divide the FOV of the sensor into $M_t^b$ parts and the births correspond to a new object being spawned in one of these FOV parts. The death hypotheses correspond to one of the $M_t^d$ objects expected to be in the FOV dying. 
Hence, a child hypothesis $H_{ij}$ of the parent $H_i$ can be an $M+1$ object hypothesis with probability $\alpha$ in exactly $M_t^b$ different ways.
The child $H_{ij}$ could have $M-1$ objects with probability $\beta$ each in $M_t^d$ different ways corresponding to the $M_t^d$ different objects dying. Thus, the child $H_{ij}$ could have $M$ objects with probability $(1-M_t^b\alpha -M_t^d\beta)$ in exactly one way (the no birth/ death case). Please see Fig. 1 for an illustration of the process.
\begin{remark}
The above development amounts to modeling the birth and death processes as independent Bernoulli processes with parameters $\alpha$ and $\beta$ respectively. Since $\alpha^2, \beta^2, \alpha\beta \approx 0$, the birth process can be modeled as just two outcomes: exactly 0 births with probability $(1-\alpha)^{M_t^b}\approx 1-M_t^b \alpha$, and exactly one birth with probability $M_t^b\alpha$. Similarly, the death process amounts to exactly 0 deaths with probability $(1-\beta)^{M_t^d} \approx 1-M_t^d\beta$, and exactly one death with probability $M_t^d \beta$. Since the two processes are independent, the joint distribution function is just a product of the two distributions and can be shown to amount to exactly 0 birth and death with probability $1-M_t^b \alpha -M_t^d\beta$, exactly one birth with probability $M_t^b\alpha$ and exactly one death with probability $M_t^d\alpha$.  The individual birth and death hypotheses then follow by noting that the births can happen in one of $M_t^b$ different ways and the deaths can happen in one of $M_t^d$ different ways.
\end{remark}

Further, the child hypothesis $H_{ij}$ can then split into further children $H_{ijk}$ where the total number of children is $A_M$, $A_{M+1}$ or $ A_{M-1}$ depending on the number of objects underlying the hypothesis $H_{ij}$, and corresponding to the various different data associations possible given the measurement $\{Z_t\}$. 
Note that the above process degenerates into the no birth and death case when $\alpha = \beta = 0$. 
Thus, we can see that the primary consequence of the birth and death process is the increase in the total number of children hypotheses. However, the equations for the multi-object filtering (with a little effort, due to the fact that the child hypotheses may have different number of objects than the parent hypothesis thereby complicating the integration underlying the prediction step) can be shown to remain unchanged. Recall Eq. \ref{FISST_M}, which is reproduced below for clarity:
\begin{align}\label{FISST_Mbd}
p(&\{X\}, (i,j)/ \mathcal{F}^t) %\nonumber\\
\hspace{1ex}=\underbrace{\frac{p(\{Z_t\}/ \{X\}, (i,j))p^-(\{X\}/ (i,j), \mathcal{F}^{t-1})}{\int p(\{Z_t\}/ \{X'\}, (i,j))p^-(\{X'\}/ (i,j), \mathcal{F}^{t-1}) d\{X'\}}}_{\mbox{updated pdf underlying} H_{ij}} \times %\nonumber\\
\frac{l_{ij} \overbrace{p_{ij}w_i}^{w_{ij}}}{\sum_{i',j'} l_{i',j'}\underbrace{p_{i'j'}w_{i'}}_{w_{i'j'}}}.
\end{align}
The only difference from the no birth and death case is, given $H_i$ is an $M-$ object hypotheses, the children hypotheses $H_{ij}$ can have $M$, $ M-1$ or $M+1$ objects underlying them, and the corresponding $p_{ij}$ value  is $1-M_t^b\alpha-M_t^d\beta$, $\beta$ or $\alpha$ respectively.
It behooves us to look  closer at the prediction equations in the birth and death case as that is the source of difference from the no birth and death case.

First, consider the case of a death hypothesis. Consider an M-object hypothesis, $H_i$, with underlying MT-pdf $\prod_k p_k(x_k)$. Suppose without loss of generality that the $M^{th}$ object dies. Then, the transition density for the multi-object system is defined as:
\begin{align}
p(\{X\}/ \{X'\}, (i,j)) = [\prod_{k=1}^{M-1} p_k(x_k/ x_k')]\delta(\phi/ x_M), 
\end{align} 
where $\delta(\phi/ x_M)$ denotes the fact that the $M^{th}$ object becomes the null object $\phi$ with probability one. Thus, the predicted MT-transition density underlying $H_{ij}$ is:
\begin{align}\label{pred_death}
p^-&(\{X\}/  (i,j), \mathcal{F}^t) =%\nonumber\\
=\int (\prod_{k=1}^{M-1} p(x_k/x_k')p(x_k'/ i, \mathcal{F}^t)) \delta(\phi/ x_M') dx_1'..dx_M'%\nonumber\\
= \prod_{k=1}^{M-1}p^-(x_k/ i, \mathcal{F}^t),
\end{align}
i.e., the predicted MT-pdf is simply the predicted pdfs of all the objects that do not die. 

Next, consider the case of a birth hypothesis $H_{ij}$ where the birthed pdf has a distribution $p_b^l(x_{M+1})$.  The transition pdf is now 
\begin{align}
p(\{X\}/ \{X'\}, (i,j)) = [\prod_{k=1}^M p_k(x_k/ x_k')]
p_{M+1}(x_{M+1}/ \phi),
\end{align}
where $p_{M+1}(x_{M+1}/\phi) = p_b^l(x_{M+1})$ denotes that the null object $\phi$ spawns an $M+1^{th}$ object with underlying pdf $p_b^l(x_M)$. It can be shown similar to above that the predicted distribution in this case is:
\begin{align} \label{pred_birth}
p^-(\{X\}/ (i,j), \mathcal{F}^t) = [\prod_{k=1}^M p_k^-(x_k/ i, \mathcal{F}^t)] p_b^l(x_{M+1}),
\end{align}
i.e., the predicted distribution of all the objects with the addition of the birth pdf $p_b^l(x_{M+1})$.

Further, each of these hypothesis split into children $H_{ijk}$ based on the possible data associations: if $H_{ij}$ is a birth hypothesis the the number of children is $A_{M+1}$, if its a death hypothesis the number of children is $A_{M-1}$ and if it is no birth or death, the number of children is $A_M$.   In particular, using the development outlined above ( where we have replaced the child notation $H_{ijk}$ by $H_{ij}$ for simplicity), we can see that the transition probability $p_{ij}$ of a child hypothesis $H_{ij}$ is:
\begin{align}\label{tp_bd}
p_{ij} = 
\begin{cases}
   \alpha p^k_{M+1},& \text{if } j\in B_{M+1, k}\\
    (1-M_t^b\alpha -M_t^d\beta){p^k_M},  & \text{if } j\in B_{M,k}\\
    {\beta}{p^k_{M-1}}, & \text{if } j\in B_{M-1, k}
\end{cases}
\end{align}
where $B_{N,l}$ represents all $N$ object hypotheses in which exactly $l$ of the objects have been associated to measurements and  
\begin{align}
p_N^l = \frac{p_D^l (1-p_D)^{N-l}}{{m \choose l} l!},
\end{align}
where $m$ is the size of the measurement return.

%\begin{center}
\begin{figure}[hbt]
\centering
\includegraphics[scale=.5]{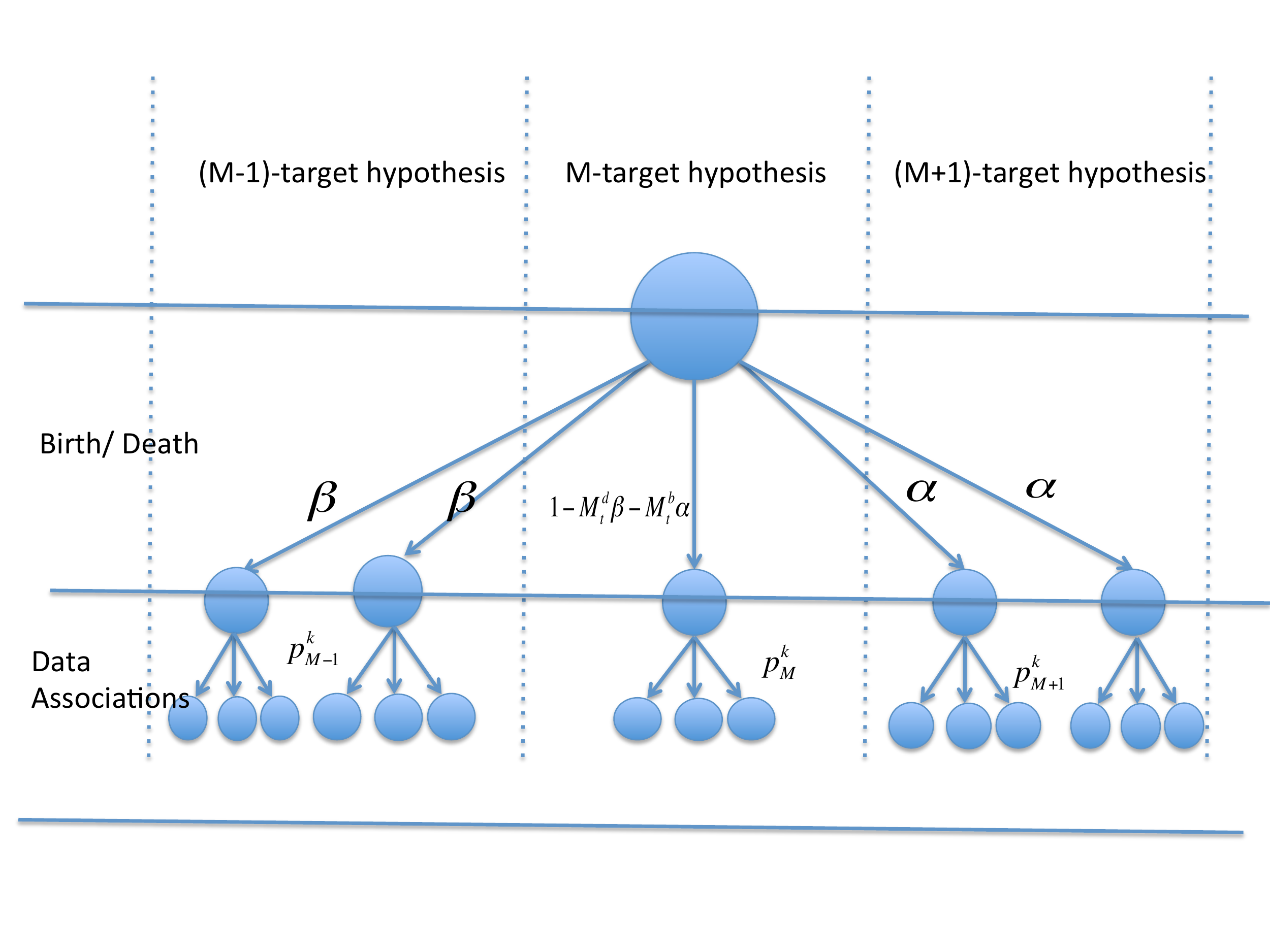}
%\label{Figure: Last}
\caption{A schematic of the splitting of the hypothesis due to birth/ death of objects and data associations. Underlying each blob is a continuous MT-pdf. }
\end{figure}
%\end{center}

The above development can be summarized as the following result:
\begin{proposition}
Given an M-object hypothesis $H_i$ and its children $H_{ij}$, the update equation for joint MT-pdf-hypothesis density function is given by Eq. \ref{FISST_Mbd}, where the only differences from the no birth or death case is that $p_{ij}$ in the equations is different according as the hypothesis $H_{ij}$ birth, death or a no birth or death hypothesis and is given by Eq.  \ref{tp_bd}, and the predicted priors required in Eq. \ref{FISST_Mbd} is calculated from Eq. \ref{pred_death} if $H_{ij}$ is a death hypothesis, Eq. \ref{pred_birth} if it is a birth hypothesis and Eq \ref{pred_nbd} if it is a no birth or death hypothesis.
\end{proposition}

\begin{remark}
The above result remains valid for more complex models of birth and death as long as the transition probabilities $p_{ij}$ can be specified. More complex models of birth and death would entail that there are significantly (combinatorially) more children hypothesis due to the births and deaths than the number here, followed by the requisite number of data association hypotheses, and would result in more complex expressions for the transition probabilities than the one in Eq. \ref{tp_bd}. As regards the relationship with MHT in the case of birth and death, the same observations as made in the no birth or death case holds: given the birth and death model, the methods result in the same exact hypotheses, except the weights of the different hypothesis are different due to the MT-likelihood normalization done in the HFISST method.
\end{remark}
 
\section{A Randomized FISST (R-FISST) Technique}
In the previous section, we have introduced the hypothesis level FISST equations and shown that they are particularly easy to comprehend and implement. However, the number of children hypothesis increase exponentially at every iteration and thus, can get computationally intractable very quickly. However, it can also be seen that most children hypotheses are very unlikely and thus, there is a need for intelligently sampling the children hypotheses such that only the highly likely hypotheses remain. In the following, we propose an MCMC based sampling scheme that allows us to choose the highly likely hypotheses.

\subsection{MCMC based Intelligent Sampling of Children Hypothesis}
Recall Eq. \ref{FISST_M'}. It is practically plausible that most children $j$ of hypothesis $H_i$ are highly unlikely, i.e., $l_{ij} \approx 0$ and thus, $w_{ij} \approx 0$. Hence, there is a need to sample the children $H_{ij}$ of hypothesis $H_i$ such that only the highly likely hypotheses are sampled, i.e., $l_{ij} >> 0$. 

\begin{remark}
Searching through the space of all possibly hypotheses quickly becomes intractable as the number of objects and measurements increase, and as time increases.
\end{remark}
\begin{remark}
We cannot sample the hypothesis naively  either, for instance, according to a uniform distribution since the highly likely hypothesis are very rare under the uniform distribution, and thus, our probability of sampling a likely hypothesis is vanishingly small under a uniform sampling distribution.
\end{remark}
Thus, we have to resort to an intelligent sampling technique, in particular, an MCMC based approach.

Given a hypothesis $H_i$, we want to sample its children according to the probabilities $\bar{p}_{ij} = w_{ij}l_{ij}$. This can be done by generating an MCMC simulation where the sampling Markov chain, after enough time has passed (the burn in period), will sample the children hypotheses according to the probabilities $\bar{p}_{ij}$.  A pseudo-code for setting up such an MCMC simulation is shown in Algorithm 1.
\begin{algorithm}
\caption{MCMC Hypothesis Sampling}
Generate child hypothesis $j_0$, set $k=0$.\\
Generate $j_{k+1} = \pi(j_k)$ where $\pi(.)$ is a symmetric proposal distribution\\
If $\bar{p}_{ij_{k+1}}> \bar{p}_{ij_k}$ then $ j_{k} := j_{k+1}; k:= k+1$;\\
else $j_k := j_{k+1}$ with probability proportional to $\frac{\bar{p}_{ij_{k+1}}}{\bar{p}_{ij_k}}$; $k = k+1$.
\end{algorithm}
In the limit, as $k\rightarrow \infty$, the sequence $\{j_k\}$ generated by the MCMC procedure above would sample the children hypotheses according to the probabilities $\bar{p}_{ij}$. Suppose that we generate $C$ highest likely distinct children hypothesis $H_{ij}$ using the MCMC procedure, then the FISST recursion Eq. \ref{FISST_M'} reduces to:
\begin{align}
w_{ij} := \frac{l_{ij}w_{ij}}{\sum_{i',j'} l_{i'j'}w_{i'j'}}, 
\end{align}
where $i'$ and $j'$ now vary from 1 to $C$ for every hypothesis $H_i$, instead of the combinatorial number $A_M$. 

Given these $M*C$ hypotheses, i.e. $C$ children of $M$ parents, we can keep a fixed number $H_{\infty}$ at every generation by either sampling the $H_{\infty}$ highest weighted hypotheses among the children, or randomly sampling $H_{\infty}$ hypotheses from all the children hypotheses according to the probabilities $w_{ij}$. 
\begin{remark}
The search for the highly likely hypotheses among a very (combinatorially) large number of options is a combinatorial search problem for which MCMC methods are particularly well suited. Thus, it is only natural that we use MCMC to search through the children hypotheses.
\end{remark}
\par
\begin{remark}
The choice of the proposal distribution $\pi(.)$ is key to the practical success of the randomized sampling scheme. Thus, an intelligent proposal choice is required for reducing the search space of the MCMC algorithm. We show such an intelligent choice for the proposal in the next section.
\end{remark}
\par
\begin{remark}
The discrete hypothesis level update Eq. \ref{FISST_M'} is key to formulating the MCMC based sampling scheme, and, hence, the computational efficiency of the R-FISST algorithm.
\end{remark} 
%%%%%%%%%%%%%%%%%%%%%%%%%%%%%%%%%%%%%%%%%%%%%%%%%%%%%%%%%%%%%%%%%%%%%%%%%%%%%%%%%%
% WES sections
\subsection{Smart Sampling Markov Chain Monte Carlo}
\begin{figure}[h]
\centering
\includegraphics[scale=.30]{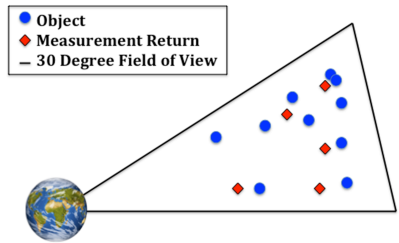}
\caption{A possible SSA event where there exists ten objects in the field of view and five measurement returns.}
\label{EXsen}
\end{figure}

\begin{figure}[h]
\centering
\includegraphics[scale=.30]{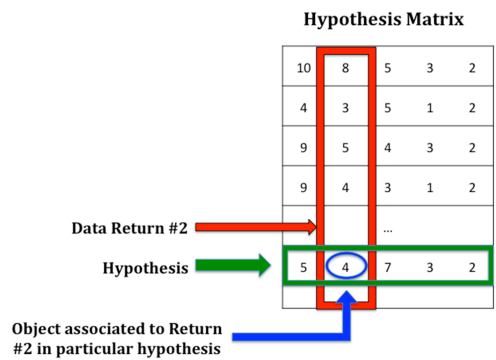}
\caption{An example of a typical Hypotheses Matrix used when exhaustively generating hypotheses. This particular matrix represents a portion of the hypothesis matrix that would be generated for the scenario in figure \ref{EXsen}.}
\label{hm}
\end{figure}

\begin{figure}[h]
\centering
\includegraphics[scale=.30]{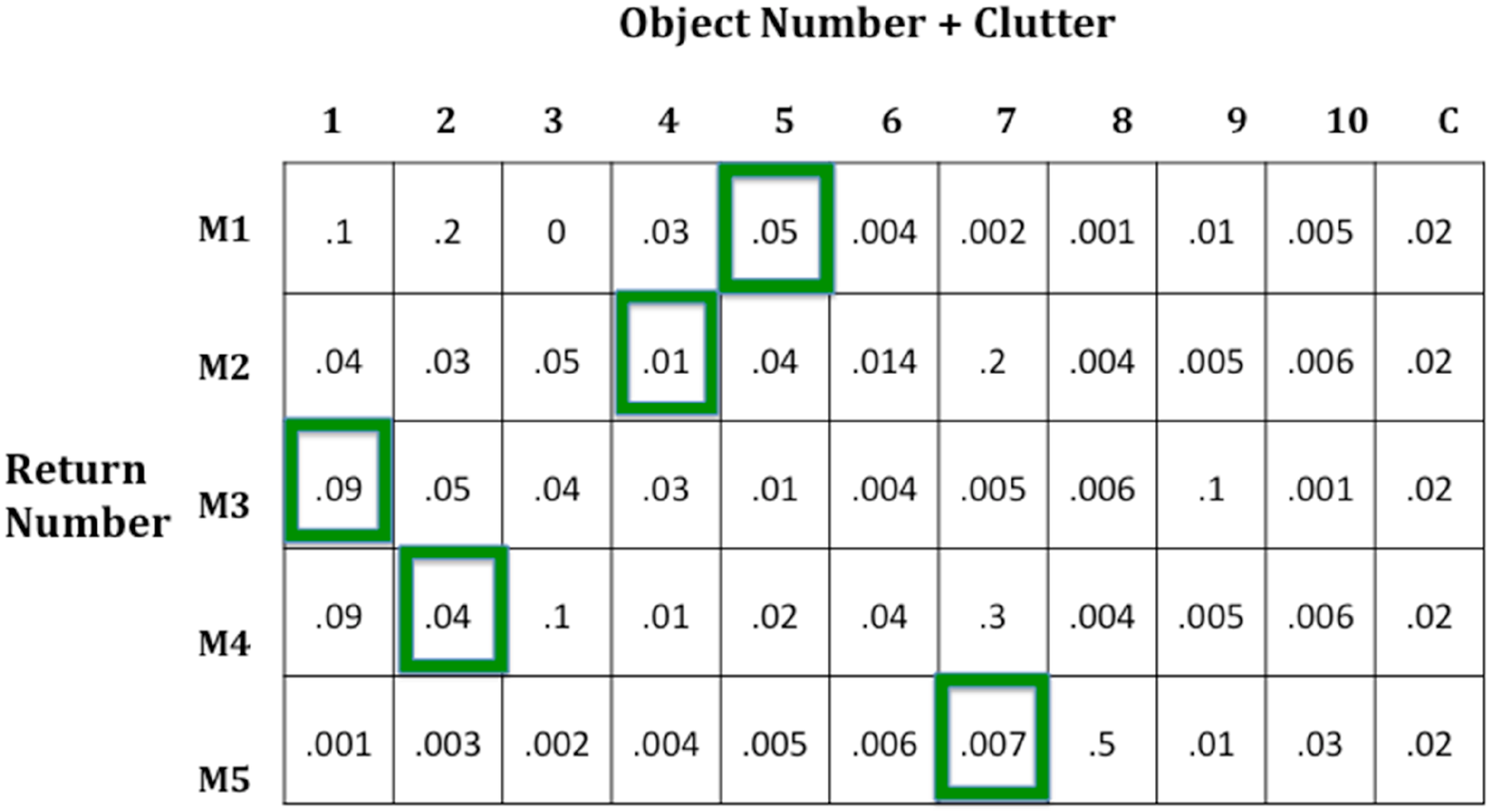}
\caption{The Data Association Matrix. Each row represents a particular measurement return. Each column is a particular association. The elements represent the likelihood of the corresponding association. Green boxes here show a visual representation of an example hypothesis.}
\label{da}
\end{figure}

In this section, we reveal the process used to perform the MCMC sampling discussed in the previous section. This process is performed at every scan to generate the highly likely children hypotheses. Consider the following SSA scenario depicted in figure \ref{EXsen}. In this scenario the belief is that there are ten objects in the field of view. The sensor then detects five measurement returns.
Typically when generating the hypotheses exhaustively one would create a matrix where each row represents a particular hypothesis. The columns of the matrix represent the measurement returns provided by the sensor. Each column entry represents the object that measurement is being associated to in the particular hypothesis. The hypothesis matrix for our example scenario would look like figure \ref{hm}. 
However, if all objects and measurements within the field of view can be associated then according to Eq. \eqref{n_da}, with $m=5$ and $M=10$, the total number of possible hypotheses would be $A_{M}=63,591$. Thus, the hypothesis matrix actually has $63,591$ rows. This illustrates the importance of a randomized approach. One can see that even with relatively low numbers of objects and measurement returns exhaustively generating and maintaining the hypotheses will cause a large computational burden. In our randomized approach we sample the rows of the hypothesis matrix based on hypothesis probability. We do this by creating a matrix we call the data association matrix, figure \ref{da}. 
The data association matrix lists the objects as the columns and the measurement return as the rows. The entries of the matrix contain the likelihood value of that particular measurement to object assignment. The last column of the matrix is dedicated to clutter and contains the likelihood that a particular measurement is associated to clutter. The dimensions of this matrix are $m\times{(M+1)}$ which is much smaller than the dimensions of the hypothesis matrix. This makes it much more practical to explore using the MCMC technique. 
\begin{remark}
The numbering of the objects and measurement returns in the data association matrix is done strictly for organization and is redone at random each time step with no record of previous numbering or labeling kept throughout scans. 
\end{remark}
\begin{remark}
Computing the data association matrix does not add any computational burden because the object to measurement likelihood is necessary in every tracking method.
\end{remark}

We start the randomized technique by creating a row vector of length $m$ containing a permutation of column numbers. The green boxes in figure \ref{da} are a visual representation of such a row vector $\left[ \begin{smallmatrix} 5&4&2&1&7 \end{smallmatrix} \right]$. This row vector is used as our first hypothesis. We then take a step in the MCMC by generating a proposed hypothesis. This is done by randomly choosing a row (Measurement) of the data association matrix to change. We then randomly sample a column (object) to associate the measurement to. If there is a conflicting assignment (i.e. a measurement is already assigned to that object) then we automatically assign the conflicting measurement to clutter. 
\begin{figure}[h]
\centering
\includegraphics[scale=.30]{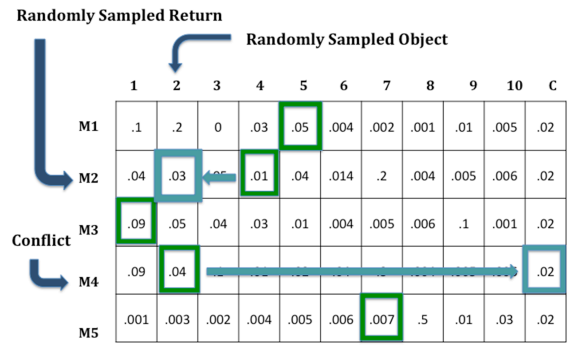}
\caption{Visualization of a single MCMC step using the Data Association Matrix. Green boxes represent the current hypothesis while blue boxes represent the changes made for the proposed hypothesis. This particular example contains a conflicting assignment with measurement return two and shows how the association is then changed to clutter.}
\label{da2}
\end{figure}
We then compare the proposed hypothesis to the current hypothesis in an MCMC fashion using a criteria which stems from the Metropolis condition 
$U[0,1]< min( 1, \frac{P_{(i,j)_{k+1}}}{P_{(i,j)_{k}}})$ where $P_{(i,j)_{k}}$ is the probability of the hypothesis at step $k$. In words, if the proposed hypothesis has a higher probability then we keep it, if not, we keep it with probability proportional to the ratio of the hypothesis probabilities. These steps are then repeated until assumed stationary distribution. We then continue walking for a user defined amount of steps and record all hypotheses sampled during these steps. The recorded hypotheses represent the highly likely hypotheses. 
\section{Applications}
This section illustrates the application of the results from the previous sections. We illustrate the R-FISST based approach to the multi-object tracking and detection problem inherent in SSA applications. In particular, we will discuss the results from a fifty-space object birth and death scenario. Our goal is to show that the aforementioned methodology allows for accurate estimation while determining the correct number of objects in an environment where the number of objects is not fixed. This will allow for the methodology to be used in both catalog update and catalog maintenance. %More in depth results including comparrissons to other tracking methods such as HOMHT are to be presented in the 
 
\subsection{R-FISST Application to a Fifty-Object Birth and Death Scenario}
In order to test the methods discussed in this paper a fifty-space object tracking and detection problem was simulated using a planar orbit model. These fifty-objects were in orbits ranging from LEO to MEO and had varying orbital properties as well as zero-mean Gaussian process noise appropriate for SSA models. The objects were simulated for at least one orbital period. That being said each object was allowed to pass completely through the field of view at least one time. The objective was to accurately track all objects given only an imperfect initial hypothesis containing some of their means and covariances. Also, to simulate a birth and death environment, the correct number of objects is initially unknown to the algorithm. In this particular example, the initial hypothesis only contains information on forty five of the fifty-objects. The five left over will be seen as objects that are randomly introduced to the environment or simply "births". This is often described as the launch of a new satellite into orbit and is not to be confused with "spawns" in which an object currently in orbit divides into two or more pieces. Spawns can be accounted for by our methodology but are not explicitly programed in this example. The R-FISST methodology must recognize the five births and provide accurate estimations of all of the objects' states and covariances. State vectors for this particular problem contain the objects' position along the $x$ and $y$ axes as well as the magnitude of their velocity in the $x$ and $y$ directions. A single noisy sensor was positioned at a fixed look direction of 15 degrees above the positive $x$-axis with a field of view of 30 degrees. The sensor was used to measure objects' position in the $x$ - $y$ plane with zero-mean Gaussian measurement noise appropriate for the application.

\begin{figure}[H]
\centering
\subfigure[Hypotheses' weights near the begining of the simulation]{
\includegraphics[width=4cm,height=2.cm]{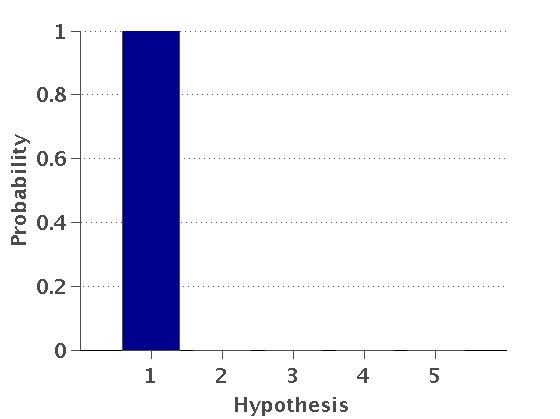}
\label{Weights1}
}
%\subfigure[Hyptheses' weights at 25 percent completion]{
%\includegraphics[width=4cm,height=2.cm]{t4}
%\label{Weights2}
%}
\subfigure[Hypotheses' weights at 50 percent completion]{
\includegraphics[width=4cm,height=2.cm]{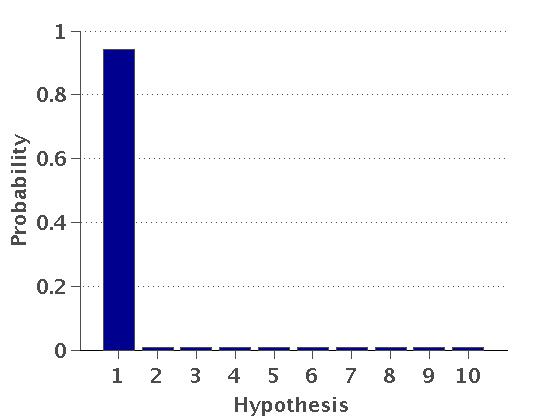}
\label{Weights3}
}
\subfigure[Hypotheses' weights at 75 percent completion]{
\includegraphics[width=4cm,height=2.cm]{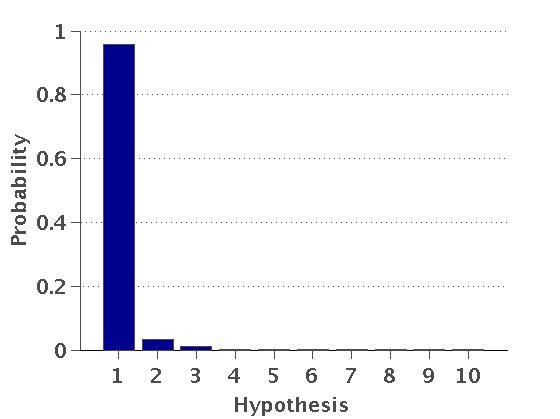}
\label{Weights4}
}
\subfigure[Hypotheses' weights at 100 percent completion]{
\includegraphics[width=4cm,height=2.cm]{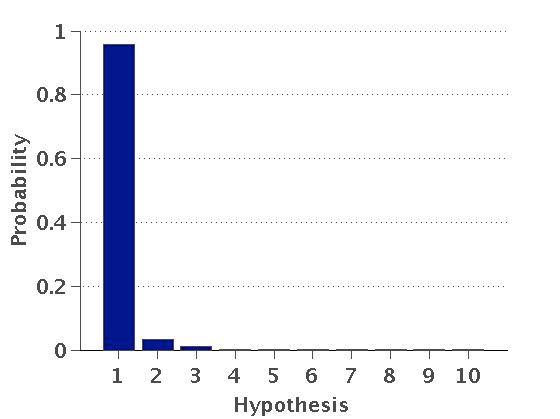}
\label{Weights5}
}
\caption{Snapshots of the hypotheses' weights throughout the simulation}
\label{Weights}
\end{figure}

\begin{figure}[H]
\centering
\subfigure[Actual Positions]{
\includegraphics[scale=.08]{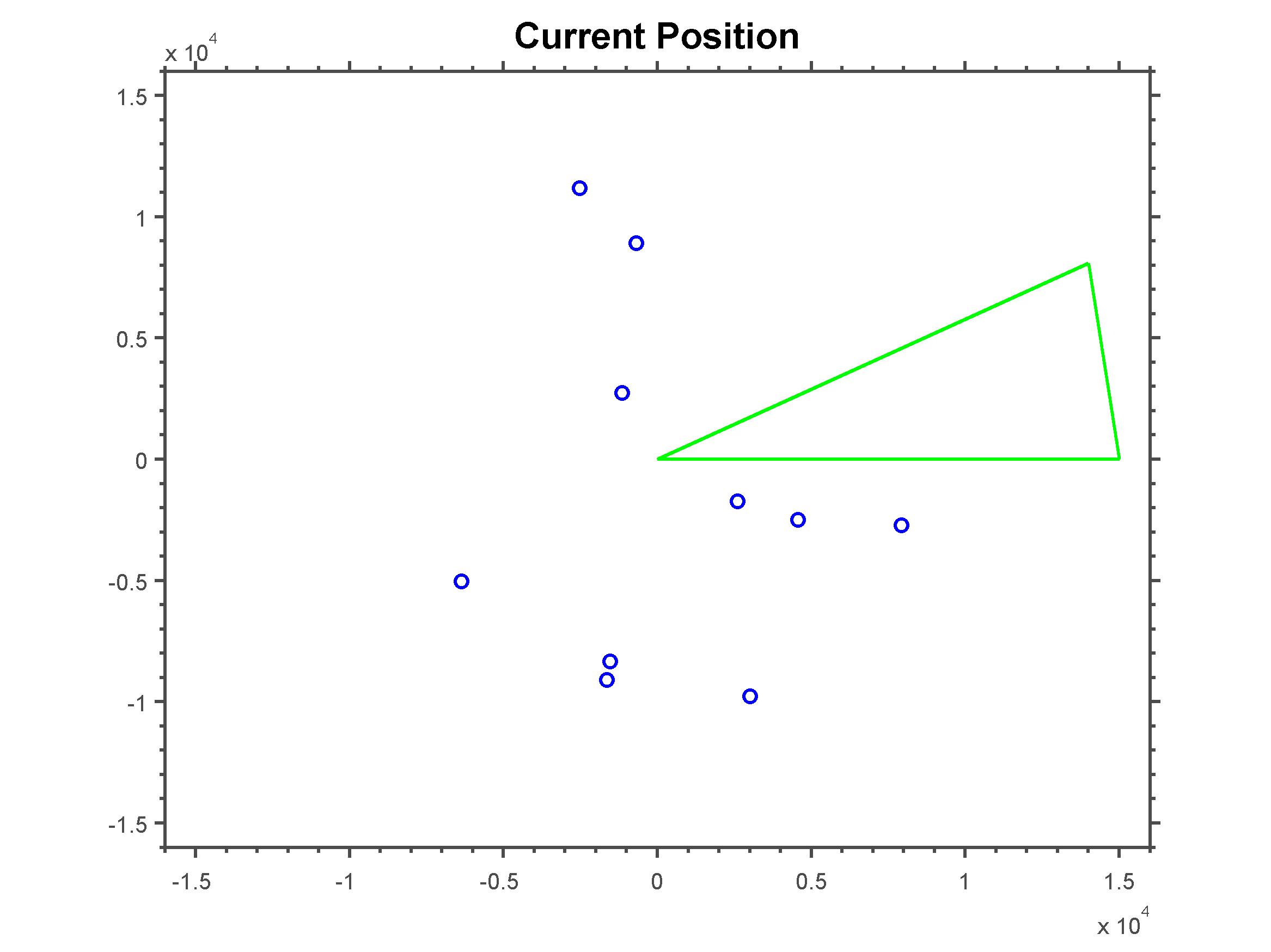}
\label{ex1}
}

\subfigure[Estimations from a correct hypothesis]{
\includegraphics[scale=.05]{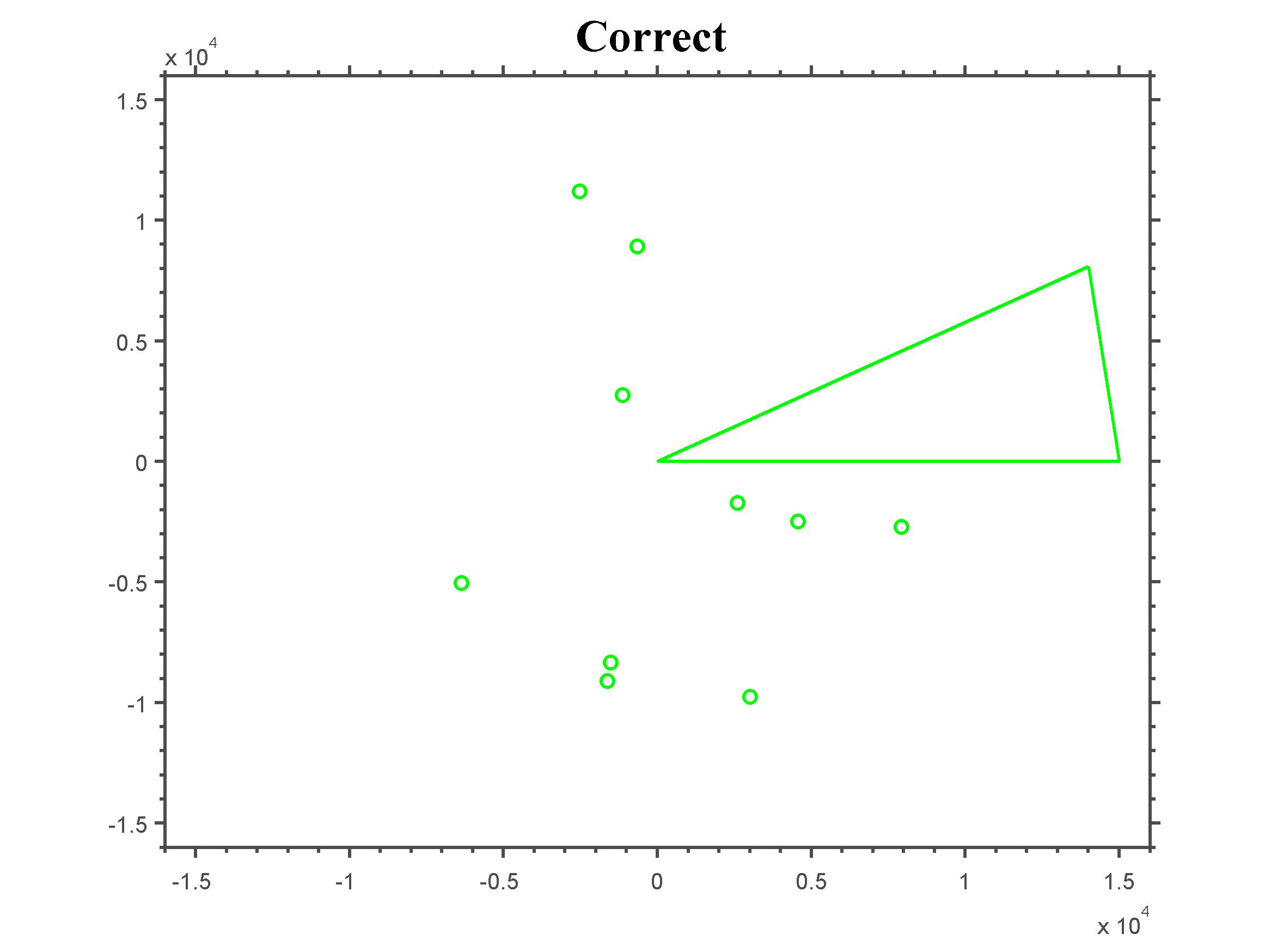}
\label{ex2}
}
\subfigure[Estimations from a hypothesis containing incorrect estimations]{
\includegraphics[scale=.05]{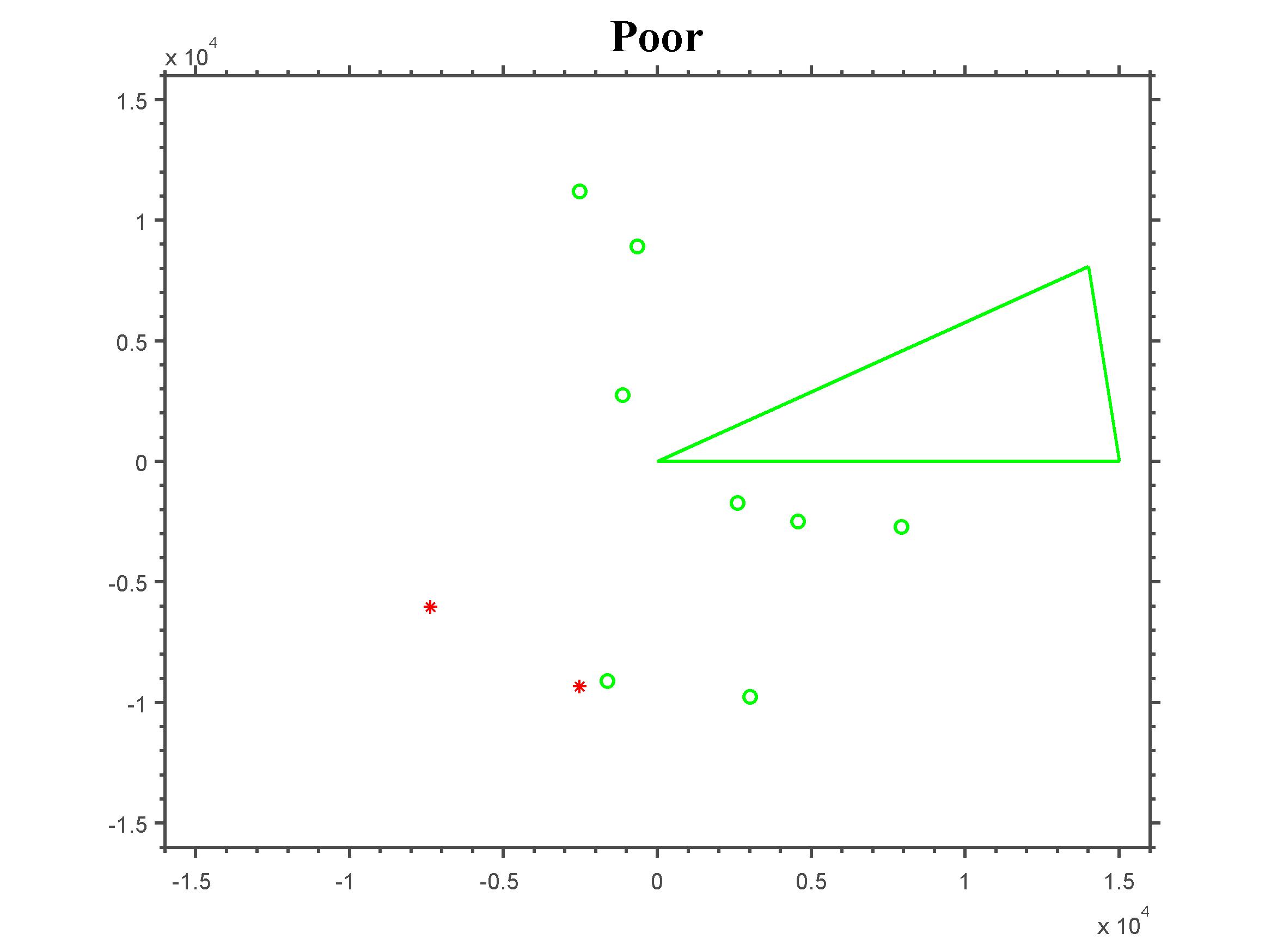}
\label{ex3}
}
\caption{An example of how estimation data is to be displayed throughout the paper. Figure \ref{ex1} shows the actual object positions while figure \ref{ex2} shows an example of a correct hypothesis and figure \ref{ex3} shows a hypothesis containing incorrect position estimations.}
\label{ex}
\end{figure}

Figure \ref{Weights} shows snapshots of the hypotheses' weights throughout the simulation. The snapshots were taken at ten, fifty, seventy-five, and one hundred percent of the total simulation time. From this figure, one can see that, in the beginning, the initial hypothesis caries all the weight. However, throughout the simulation the number of maintained hypotheses (shown on the x-axis of the graphs in Figure \ref{Weights}) varies as does the weights of those hypotheses. The number of hypotheses maintained has a direct correlation to the number of recent ambiguities in the field of view. Ambiguities occur when one or more measurement returns can be associated to multiple objects in the field of view.
Before discussing the state estimations, figure \ref{ex} shows an example of how the estimation data is to be presented. The figure shows the actual positions of the objects labeled "Current Position" and the estimated positions from two seperate hypotheses. If the estimated position for a particular object is within an error bound then a green circle will represent the object’s position otherwise a red star will represent the object’s position. %The first hypothesis labeled "Correct" shows a situation where the hypothesis corresponded to estimates of all object positions within the error bound. The second hypothesis labeled "Poor" shows a situation where the hypothesis corresponded to some estimates that were accurate and others that were above the error tolerance. 

\begin{figure}[h]
\centering
\subfigure[Actual Object Positions]{
\includegraphics[scale=.05]{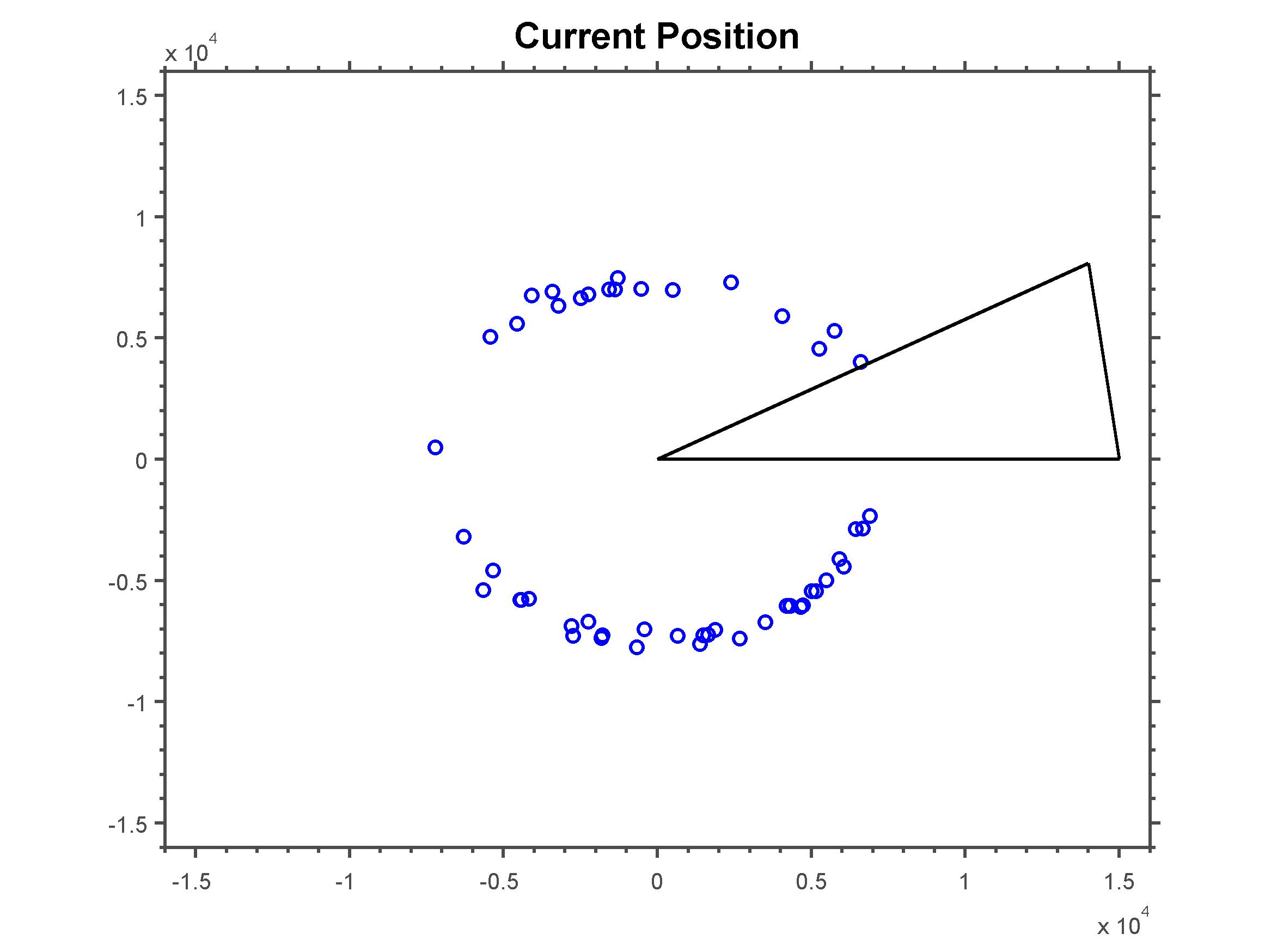}
\label{SE1}
}
\subfigure[Estimation from the top hypothesis at 10 percent completion]{
\includegraphics[scale=.05]{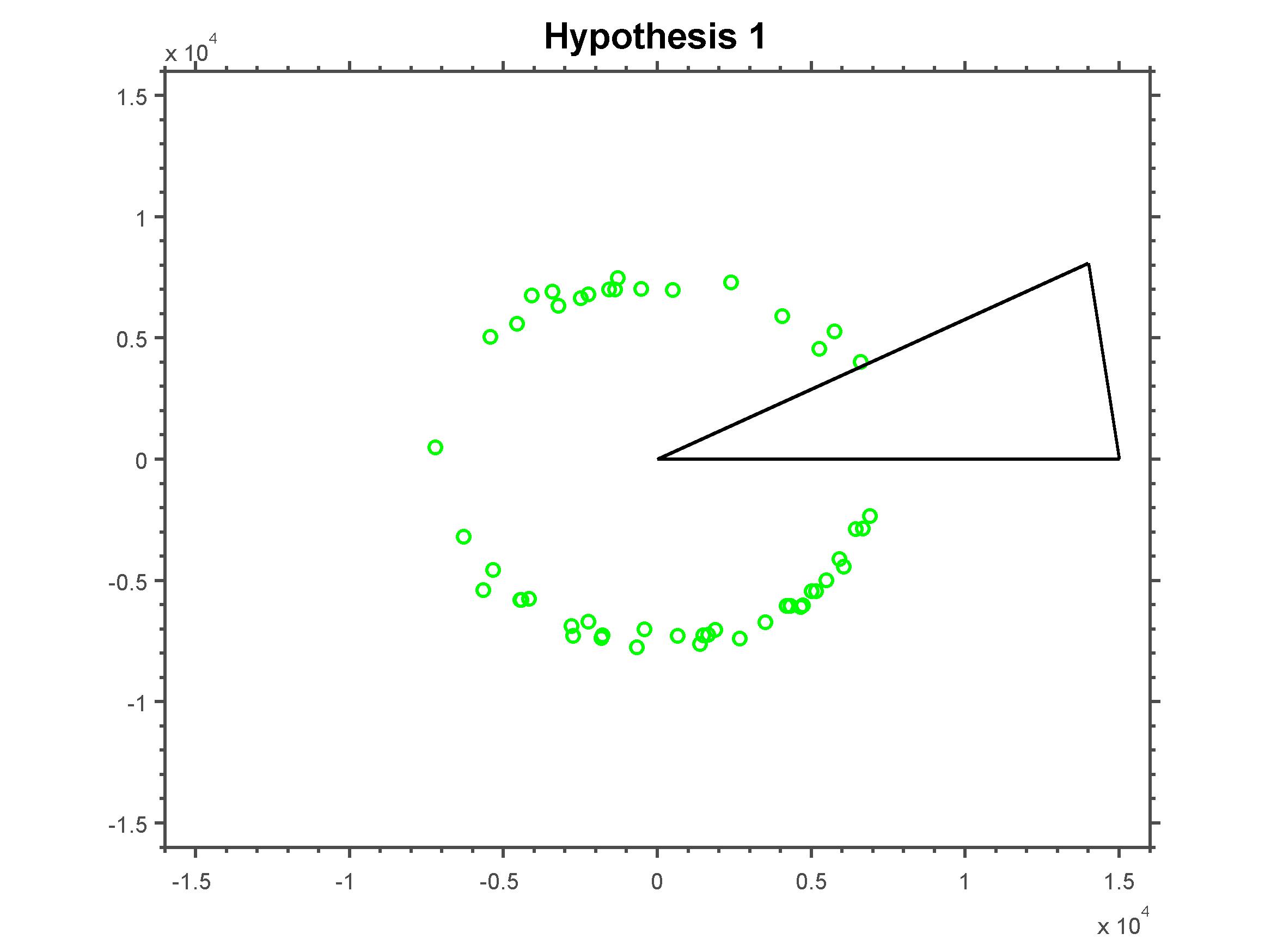}
\label{SE1B}
}
\caption{Snapshots of the actual states (black) and the estimated states from the top hypotheses (green) at 10 percent completion. Axes in tens of thousands of kilometers}
\label{TrialEst1a}
\end{figure}

\begin{figure}[h]
\centering
\subfigure[Actual Object Positions]{
\includegraphics[scale=.05]{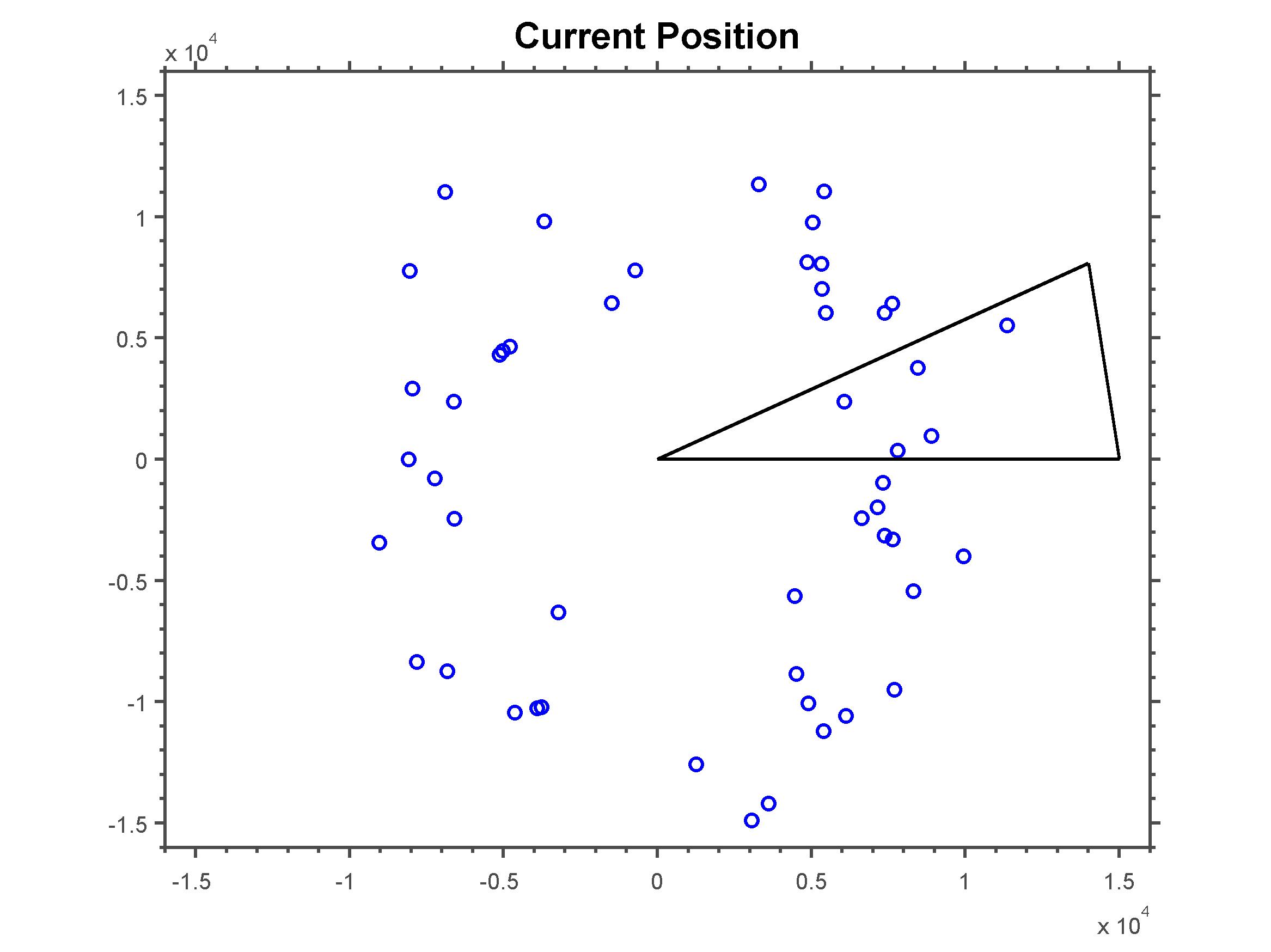}
\label{SE2}
}
\subfigure[Estimation from the top hypothesis at 50 percent completion]{
\includegraphics[scale=.05]{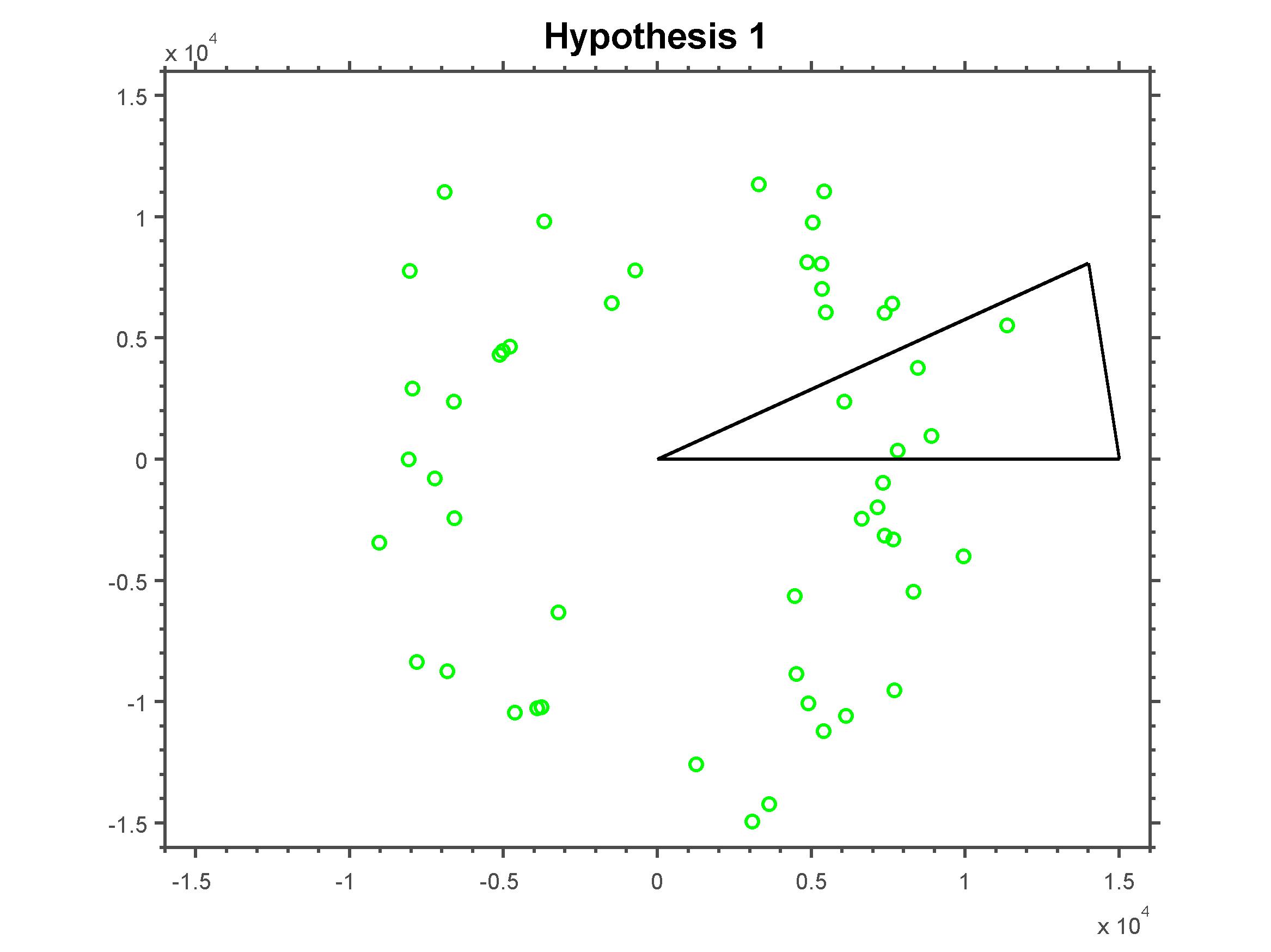}
\label{SE2b}
}
\caption{Snapshots of the actual states (black) and the estimated states from the top hypotheses (green) at 50 percent completion. Axes in tens of thousands of kilometers}
\label{TrialEst1b}
\end{figure}

\begin{figure}[h]
\centering
\subfigure[Actual Object Positions]{
\includegraphics[scale=.05]{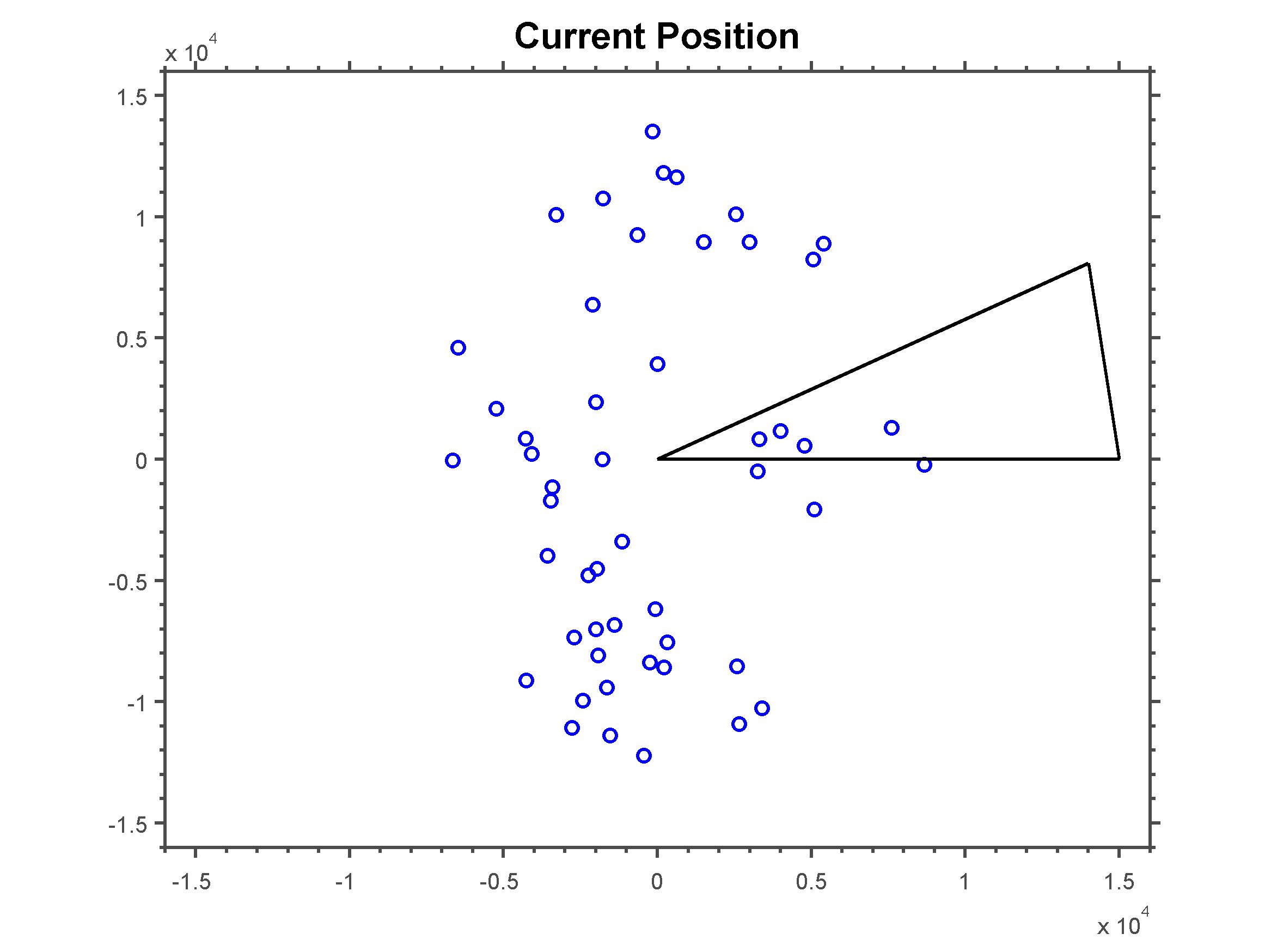}
\label{SE3}
}
\subfigure[Estimation from the top hypothesis at 100 percent completion]{
\includegraphics[scale=.05]{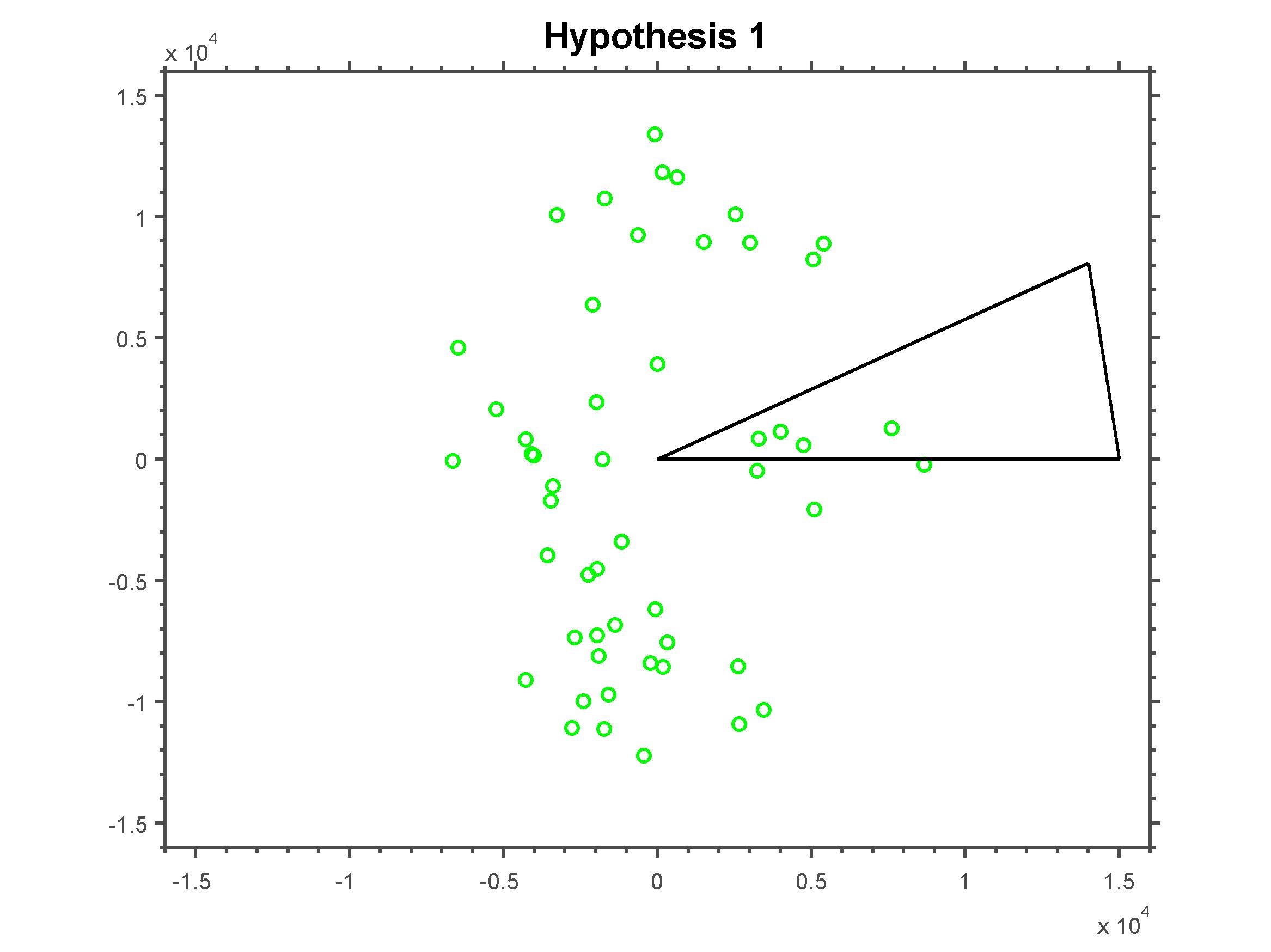}
\label{SE3b}
}
\caption{Snapshots of the actual states (black) and the estimated states from the top hypotheses (green) at 100 percent completion. Axes in tens of thousands of kilometers}
\label{TrialEst1c}
\end{figure}

Figure \ref{TrialEst1a}-\ref{TrialEst1c}, the snapshots show the actual positions (blue) against the estimated position from the top hypotheses (green). The black lines bound the field of view. These snapshots were taken at the same time intervals as in figure \ref{Weights} and thus the estimates throughout figures \ref{TrialEst1a}-\ref{TrialEst1c} are taken from the hypotheses with the highest weights in figure \ref{Weights}. Notice in figures \ref{TrialEst1a}-\ref{TrialEst1c} there are no instances of red stars. This is particularly important because it shows that the hypotheses accurately estimated object positions throughout the simulation. Hence, the R-FISST approach accurately tracked and detected the fifty-objects.

\begin{figure}[]
\centering
\includegraphics[scale=.36]{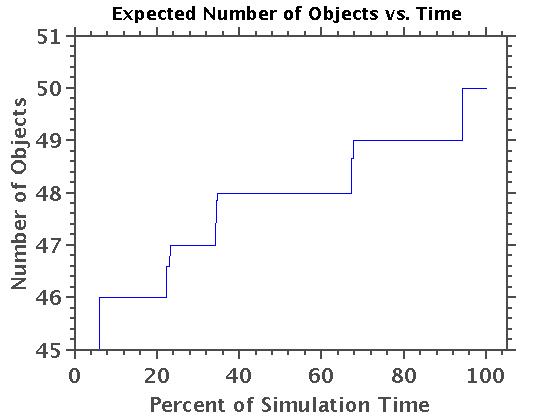}
\caption{The expected number of objects throughout the simulation. The graph shows intermediate values during the transition periods. Over time it can be seen that weight shifts from the forty-five object hypotheses to the fifty-Object hypotheses}
\label{Ws}
\end{figure}

Lastly, figure \ref{Ws} shows a visual representation of how weight shifts from the forty-five object assumption to the fifty-object assumption. The $x$-axis represents the simulation time in percent completed. The $y$-axis represents the expected number of objects. The expected number of objects is found by summing the weights of all hypotheses containing the same number of objects. The magnitudes of these summations are then compared to determine the expected number of objects. It is important to note that weight seems to be handed off in a single file fashion until the fifty-object assumption accumulates all the weight toward the end of the simulation. 

\section{Comparisons}
In this section, we show a comparison between our randomized method RFISST and a well-known tracking method called HOMHT. This comparison helps us illustrate two main points. The first being that the accuracy of the estimation provided by the RFISST method is either equal to or better than that of HOMHT but never worse. We will achieve this by showing side-by-side estimations from both methods. The second point is seen when the number of hypotheses rapidly escalates due to a large number of objects and/or a large number of measurement returns. Such occurrences happen often in SSA and for many reasons, for example, when a debris field crosses the sensor's field of view. In these situations HOMHT fails because it is computationally impossible to generate such a large number of hypotheses. Due to the randomized scheme, the RFISST methodology continues to perform in such scenarios.

\subsection{Comparison between RFISST and HOMHT: SSA Tracking}
In order to compare both methods we simulated a fifteen-space object tracking and detection problem. Each object was given a random planar orbit ranging between LEO and MEO with unique orbital properties and zero-mean Gaussian process noise appropriate for SSA models. The objects were simulated for long enough to where the object with the largest period would be able to complete at least one orbit. Thus each object was allowed to pass completely through the field of view at least one time. In this particular simulation we initialized all orbits to begin within the field of view. In order to achieve an apples to apples comparison we used this simulation to test both methods. The goal of each method would be to accurately track each object given only an imperfect initial hypothesis containing the objects' mean and covariance as well as measurement returns from a single noisy sensor. State vectors for this problem consisted of the objects' position along the $x$ and $y$ axes as well as the magnitude of their velocity in the $x$ and $y$ directions. The single noisy sensor was positioned at a fixed look direction of 15 degrees above the positive $x$-axis with a field of view of 30 degrees. The sensor was used to measure objects' position in the $x$ - $y$ plane with zero-mean Gaussian measurement noise appropriate for the application. An Extended Kalman Filter (EKF) was used in conjunction with each method to compute the underlying state and covariance updates. That being said both methods will produce the same estimation given the correct hypotheses were generated throughout the simulation. 
\begin{figure}[]
\centering
\subfigure[]{
\includegraphics[scale=.08]{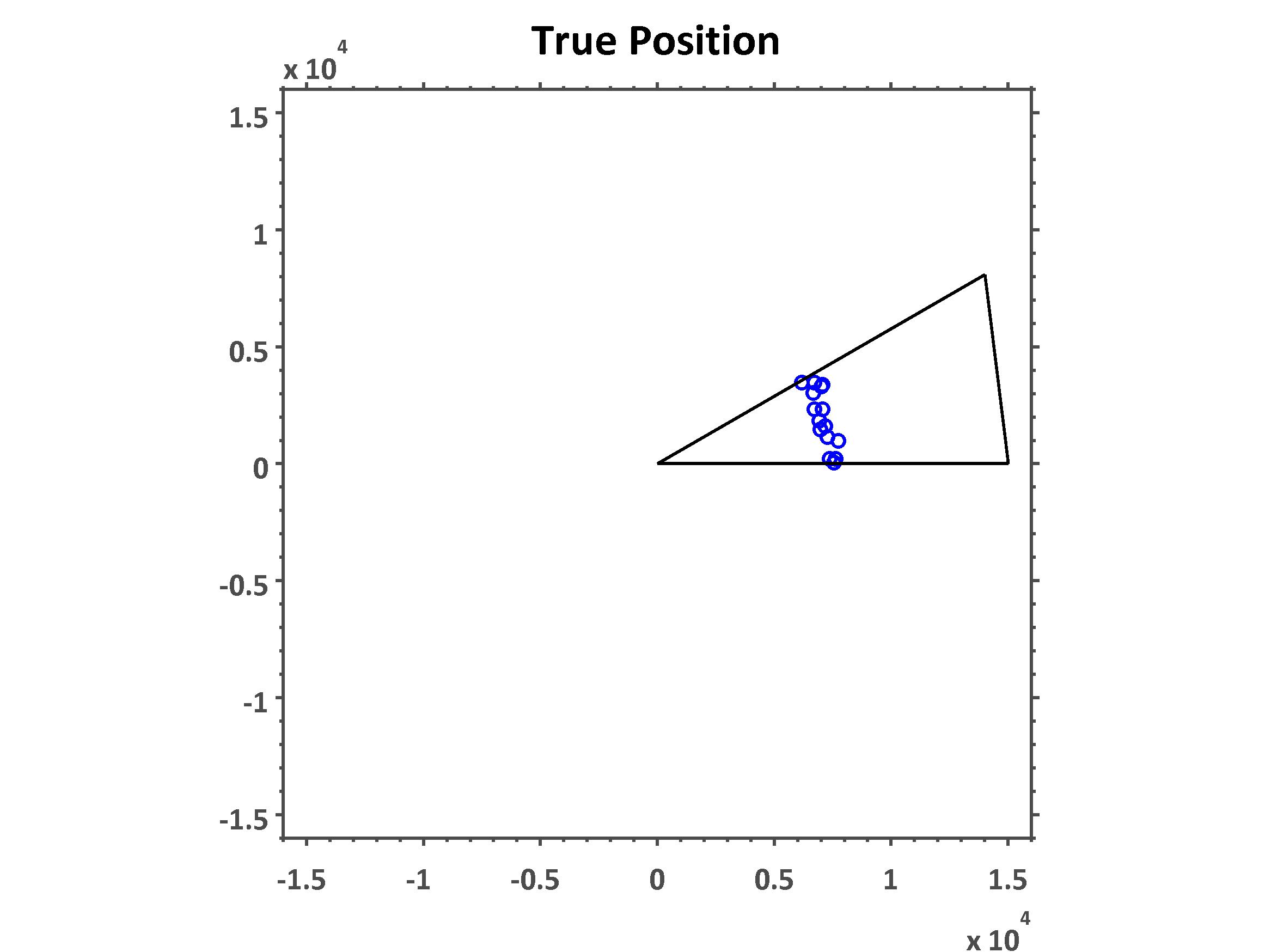}
\label{c1a}
}

\subfigure[]{
\includegraphics[scale=.05]{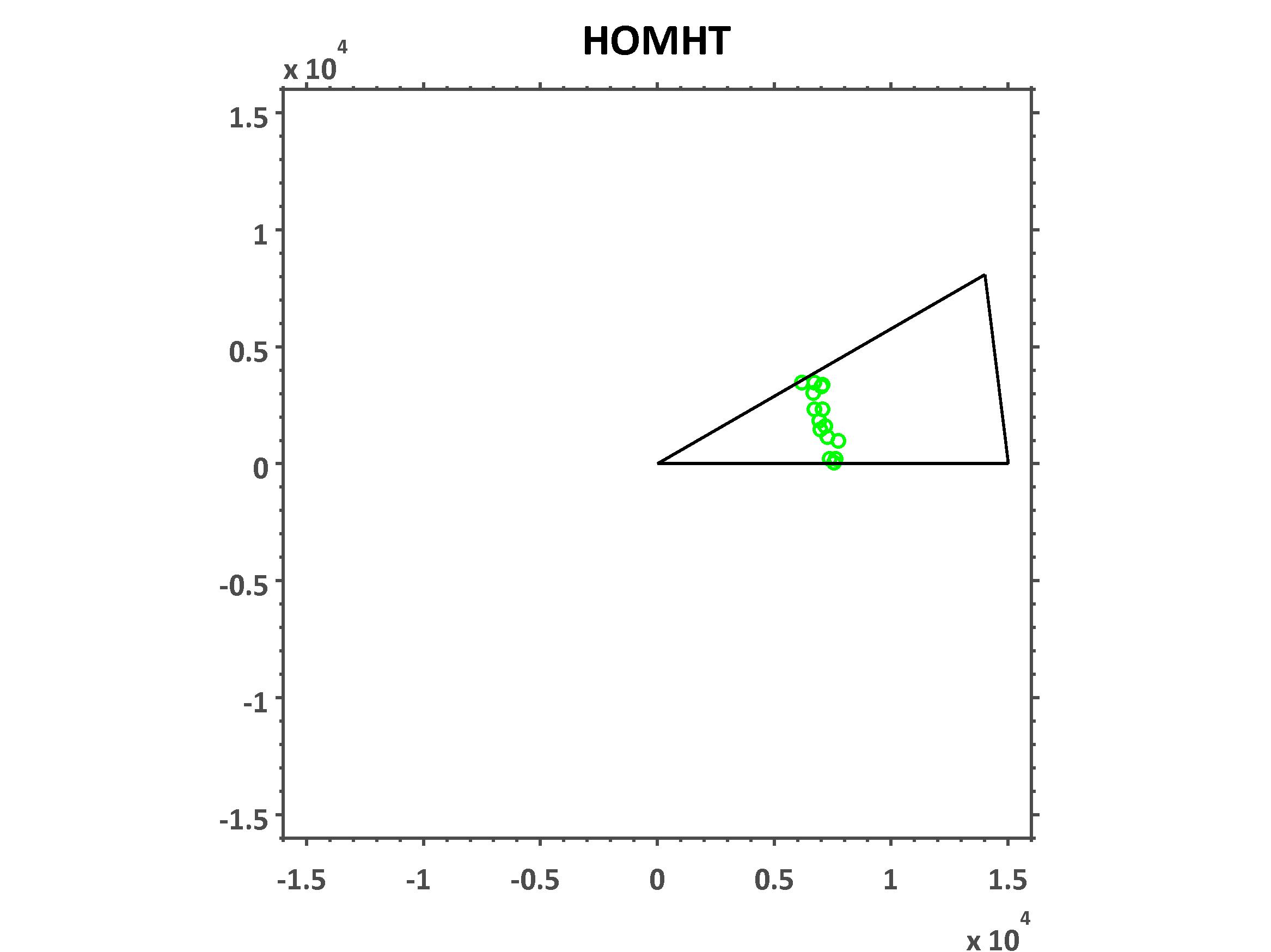}
\label{c1b}
}
\subfigure[]{
\includegraphics[scale=.05]{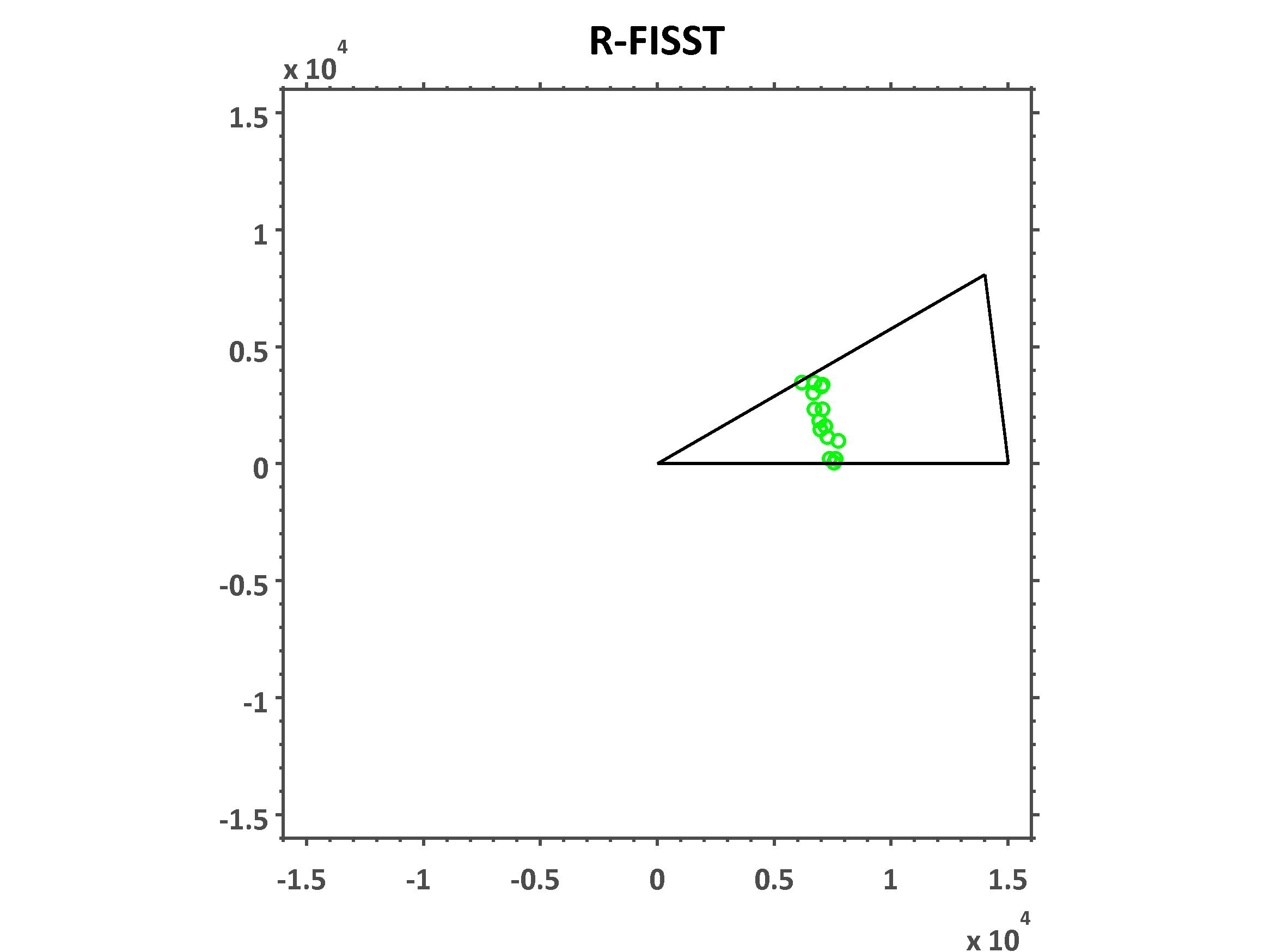}
\label{c1c}
}
\caption{Estimation at the beginning of the simulation}
\label{c1}
\end{figure}

\begin{figure}[]
\centering
\subfigure[]{
\includegraphics[scale=.08]{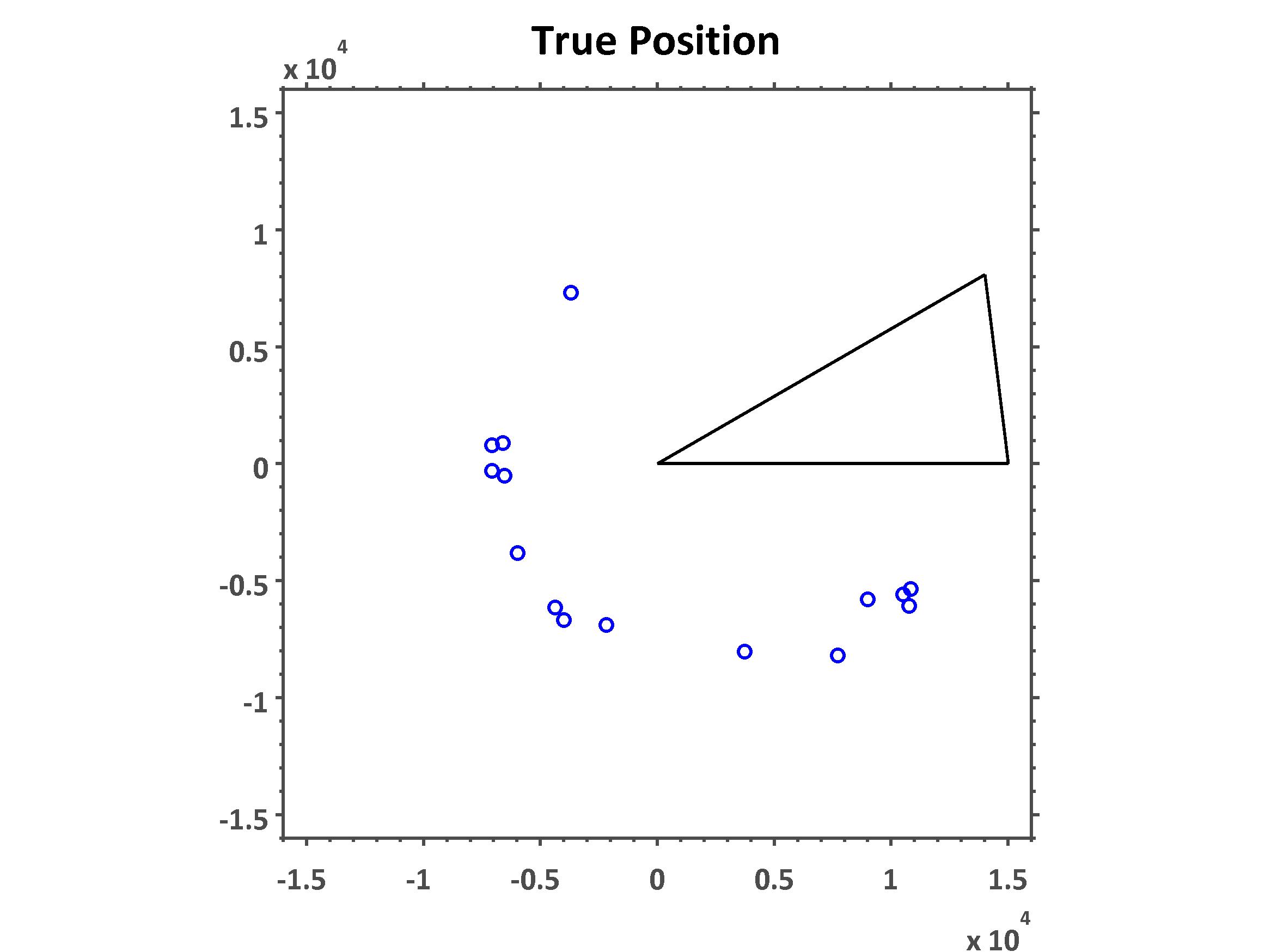}
\label{c2a}
}

\subfigure[]{
\includegraphics[scale=.05]{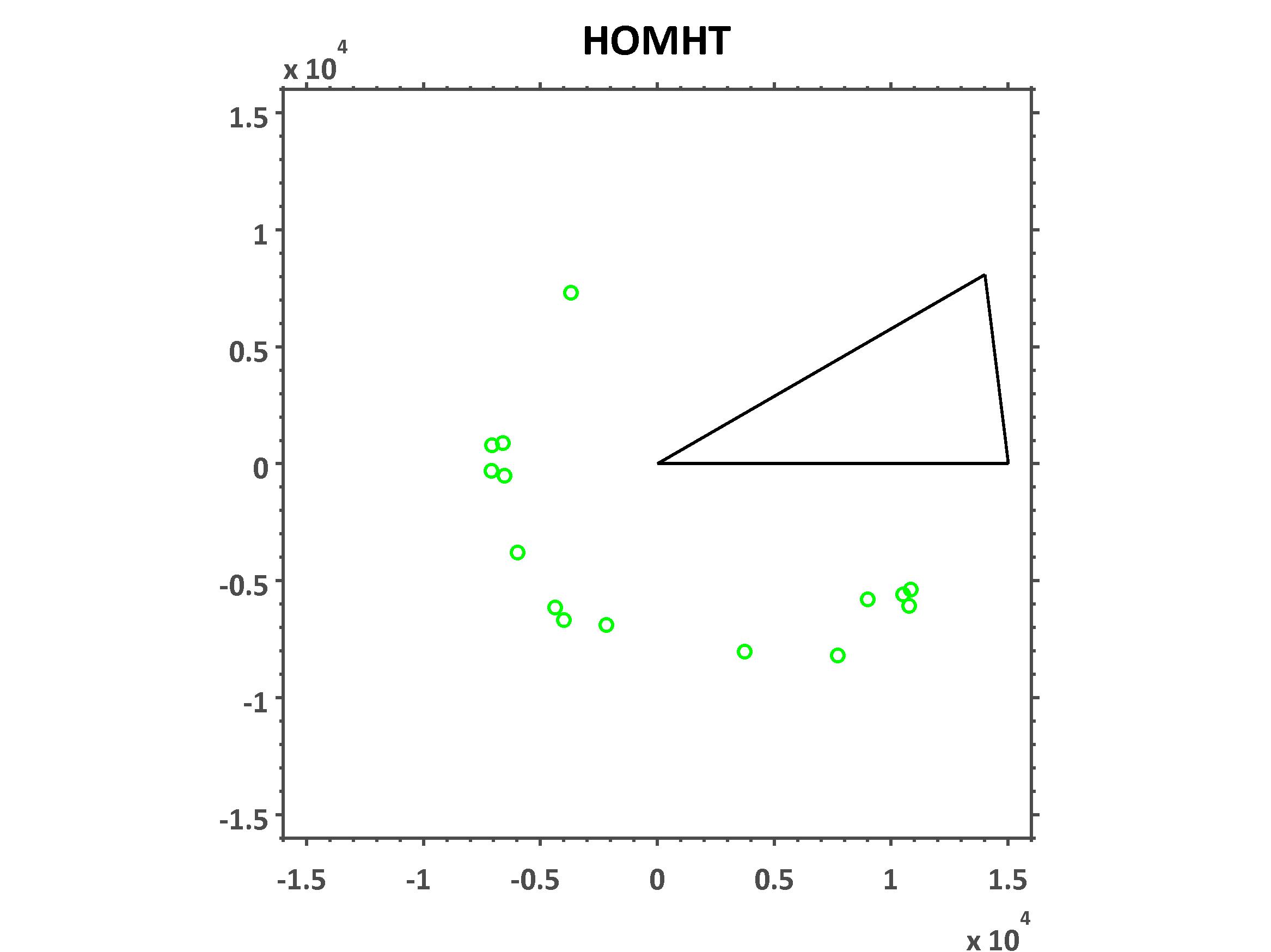}
\label{c2b}
}
\subfigure[]{
\includegraphics[scale=.05]{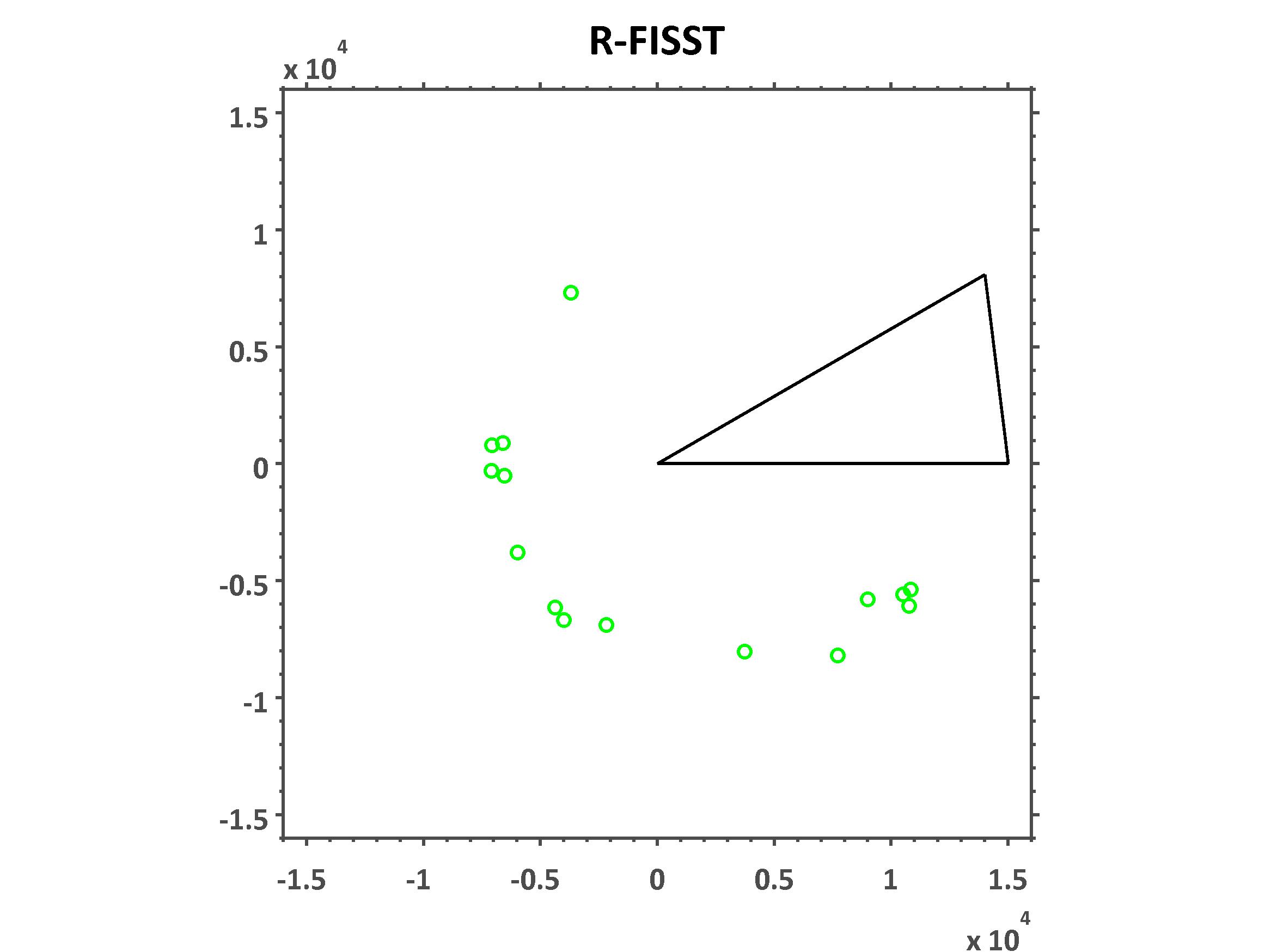}
\label{c2c}
}
\caption{Estimation at 50 percent completion of the simulation}
\label{c2}
\end{figure}

\begin{figure}[]
\centering
\subfigure[]{
\includegraphics[scale=.08]{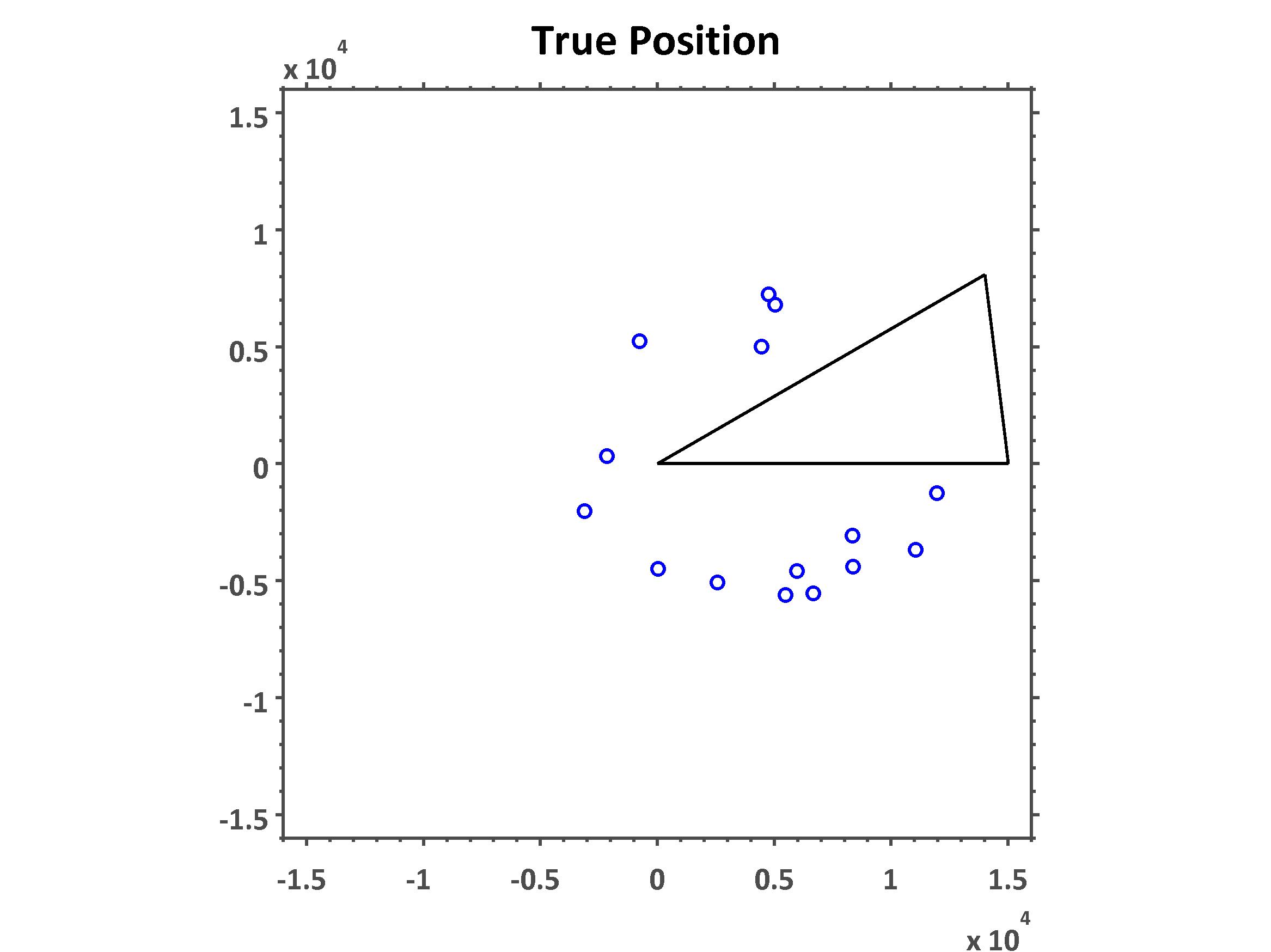}
\label{c3a}
}

\subfigure[]{
\includegraphics[scale=.05]{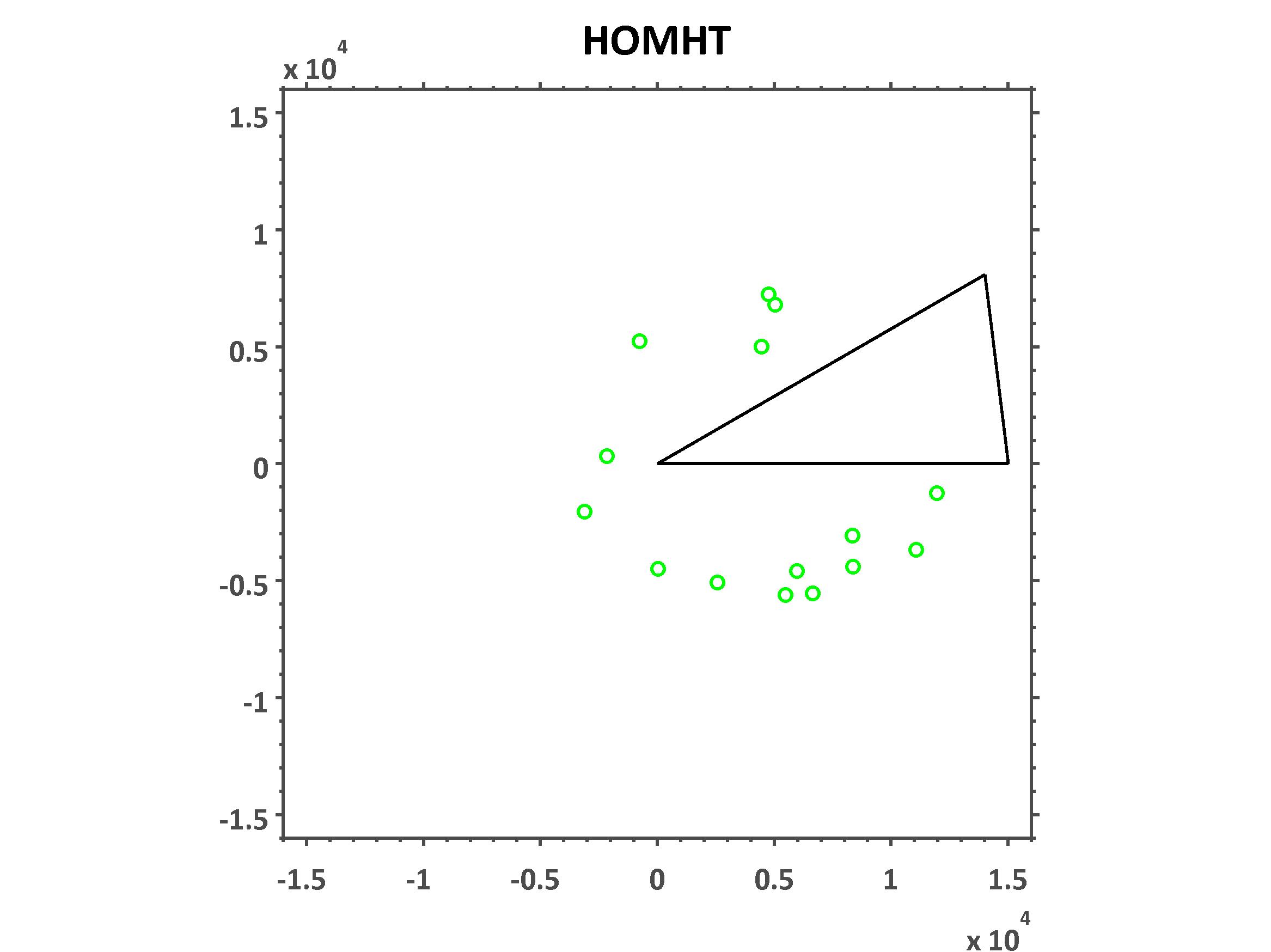}
\label{c3b}
}
\subfigure[]{
\includegraphics[scale=.05]{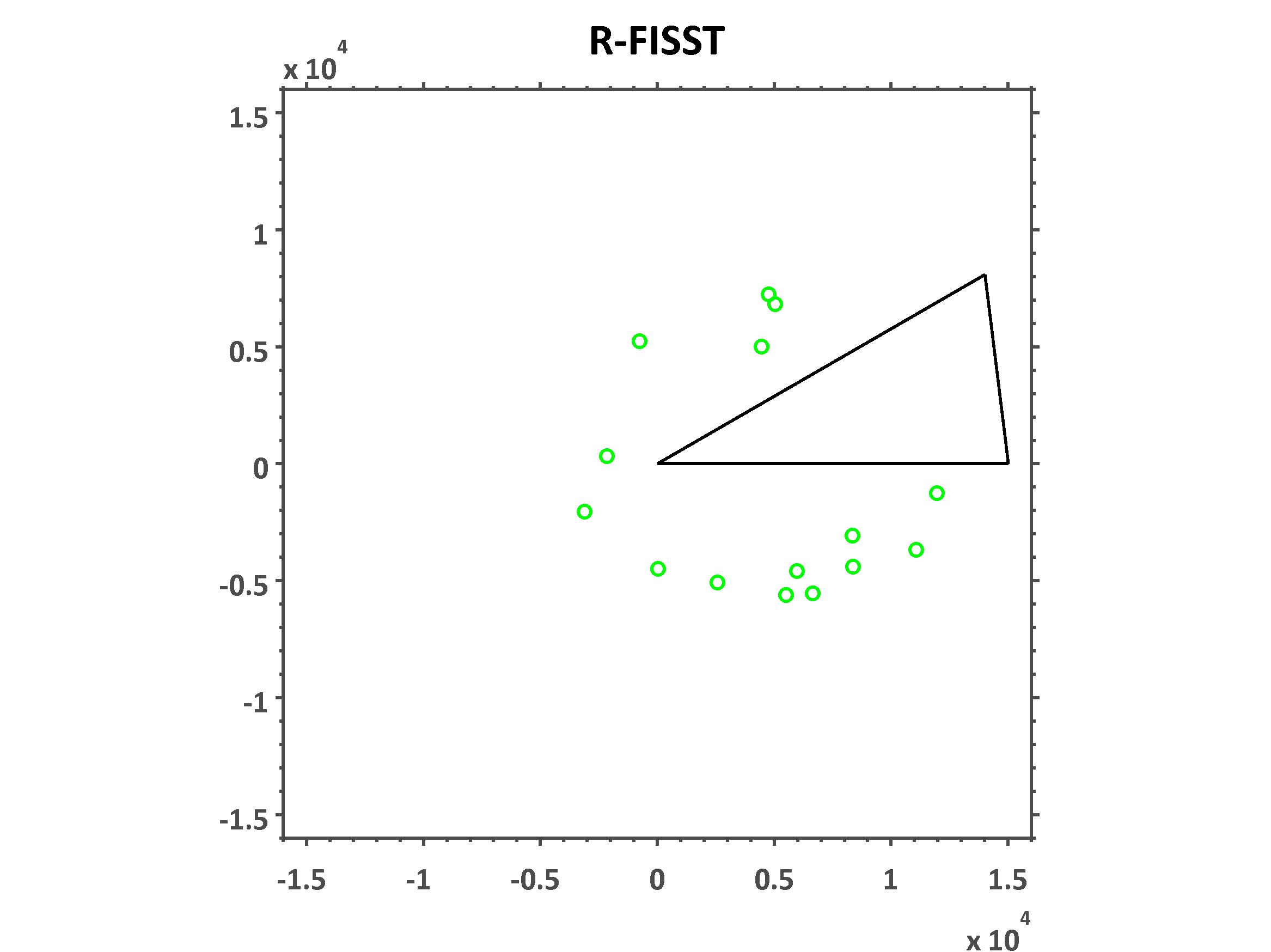}
\label{c3c}
}
\caption{Estimation at the end of the simulation}
\label{c3}
\end{figure}

Figures \ref{c1}-\ref{c3} are snapshots of the simulation at the beginning middle and end of the simulation time. Each snapshot shows the true positions of the objects and the position estimations provided by both the HOMHT method and the RFISST method. These figures are provided to illustrate that the methods accurately track the objects. Furthermore, it shows that each method maintained the correct hypothesis throughout the simulation. If either method was unable to generate the correct hypothesis then the position estimations would be incorrect. These incorrect position estimates would be seen as red stars in the snapshots. 
\begin{figure}[H]
\centering
\includegraphics[scale=.3]{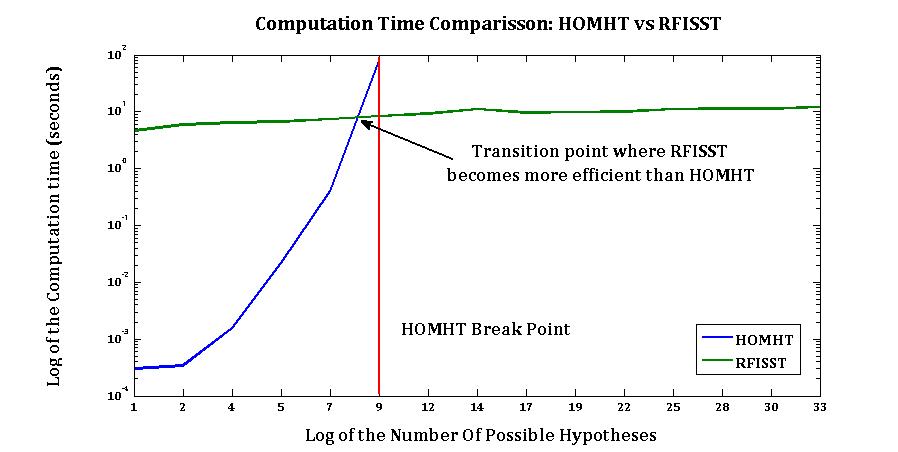}
\caption{The hypotheses generation time for both methods as shown on a log y-axis scale. }
\label{Logy}
\end{figure}
%Discuss computation time plots

\subsection{Comparison between RFISST and HOMHT: Computation Time}
When tracking large numbers of objects the majority of the computational burden lies in the hypothesis generation. HOMHT uses an exhaustive approach to generate the hypotheses. This exhaustive approach has its pros and cons. For example, generating all the hypotheses guarantees that the hypothesis containing the correct data associations is sampled. Also this exhaustive approach is very easy to implement. On the other hand, as the number of hypotheses grows so does the burden placed on generating them. This can be seen as an increase of computation time. Even with proper gating and pruning methods the number of hypotheses can be so large that generating them would exceed the computers memory heap space making it computationally intractable. It is at this point that we say the HOMHT method breaks. Using our randomized approach we never generate all hypotheses, which allows us to handle scenarios with very high number of possible hypotheses. However this method is more difficult to implement and must be tuned to guarantee that the hypothesis containing the correct data association is sampled. Figure \ref{Logy} shows a comparison of the computation times for hypothesis generation of HOMHT and RFISST. The y-axis is the log of the computation time in seconds while the x-axis is the log of the total number of possible hypotheses being generated. HOMHT is represented by the blue line and resembles an exponential curve. RFISST is represented by the green line and resembles a linear growth. It can be seen that at first for low numbers of possible hypotheses HOMHT performs faster. However, as the number of possible hypotheses grows into the tens of thousands the RFISST method becomes more efficient. Furthermore, as the number of possible hypotheses grows into the hundred millions HOMHT struggles to generate the hypotheses and eventually breaks. RFISST can generate the correct hypotheses even as the number of possible hypotheses grows to the order of $10^{33}$. In each of these simulations the RFISST MCMC methodology was able to sample the correct hypotheses using only $100,000$ steps in the MCMC. 

\section{Conclusion}
In this paper, we have presented an alternate hypothesis based derivation of the FISST recursions for multi-object tracking. We have also introduced a randomized version of the FISST recursion, called R-FISST, which provides a computationally efficient randomized solution to the full FISST recursions via an MCMC sampling of likely children hypotheses. We have also proposed a unification of the hitherto deemed different FISST and MHT methodologies for multi-target tracking. We have shown the capability of the R-FISST method using a fifty-object birth and death SSA scenario. We  showed that given the same underlying tracking filter, the HOMHT and RFISST produce the similar tracking performance. We further showed that in situations involving high numbers of possible hypotheses, the exhaustive generation of hypotheses typically used in HOMHT becomes computationally intractable while the RFISST method continues to perform well. Currently we are looking to apply our method to real data and develop large scale GPU based implementation that can scale to realistic scenarios. We also intend to look at the integration of sensor tasking into the tracking methodology such that the ambiguities inherent in the problem can be minimized.

%%%%%%%%%%%%%%%%%%%%%%%%%%%%%%%%%%%%%%%%%%%%%%%%%%%%%%%%%%%%%%%%%%%%%%%%%%%%%%%%
%%%%%%%%%%%%%%%%%%%%%%%%%%%%%%%%%%%%%%%%%%%%%%%%%%%%%%%%%%%%%%%%%%%%%%%%%%%%%%%%%%
\section{ACKNOWLEDGMENTS}
\label{.7}
This work is funded by AFOSR grant number: FA9550-13-1-0074 under the Dynamic Data Driven Application Systems (DDDAS) program.
%\input{reference}
%%%%%%%%%%%%%%%%%%%%%%%%%%%%%%%%%%%%%%%%%%%%%%%%%%%%%%%%%%%%%%%%%%%%%%%%%%%%%%%%%%
\appendix
In this appendix, we shall show the equivalence of the FISST equations and the hypothesis level derivation of this paper. Let the FISST pdf initially of the form:
\begin{align}
p(\{X\}) = \sum_{i=1}^N \omega_i p_i(\{X\}),
\end{align}
where for each $i$, the pdf $p_i(\{X\})$ is a Multi-Target FISST pdf of the form:
\begin{align}
p_i(\{X\}) = \sum_{\bar{\sigma}} p^i_{\sigma_1}(x_1)p^i_{\sigma_2}(x_2)\cdots p^i_{\sigma_n} (x_n),
\end{align}
where $\bar{\sigma} = \{\sigma_1, \sigma_2 \cdots \sigma_n\}$ denotes all possible permutations of the indices $\{1, 2 \cdots n\}$, and $\omega_i$ is a non negative number such that $\sum_i \omega_i = 1$. \\

In the following, for notational simplicity and clarity, we assume that the MT-pdf corresponding to the $i^{th}$ component is a 2-target pdf. The work essentially carries over directly to the n-target case at the expense of more notation. Thus, the MT-pdf is given by:
\begin{align}
p_i(\{x_1, x_2\}) = p^i_1(x_1)p^i_2(x_2) + p^i_1(x_2)p^i_2(x_1).
\end{align}
In the case of the hypothesis level derivation, call it H-FISST,  $\omega_i$ corresponds to the wt of the $i^{th}$ hypothesis and the underlying MT-pdf is given by $p^i_1(x_1)p^i_2(x_2)$, i.e, the H-FISST MT-pdf is the same as the FISST MT-pdf sans the permutation of the arguments $(x_1, x_2)$.\\

Next, we look at the prediction step of the FISST equations. We will assume no target birth or death but the situation is very similar even in the case of birth and death, and can be derived analogous to the following.\\
The MT-transition FISST pdf is given by:
\begin{align}
p(\{x_1,x_2\}/ \{x_1', x_2'\}) = p(x_1/x_1')p(x_2/x_2') + p(x_1/x_2')p(x_2/ x_1'),
\end{align}
where $p(x/x')$ is the transition pdf of a single target which assume here to be the same for all targets. Again, it may be extended to different classes of transition pdfs at the expense of more notation. The FISST predicted pdf for the $i^{th}$ component is then given by:
\begin{align}
p^-_i(\{x_1, x_2\}) = \nonumber\\
\frac{1}{2!} \int (p(x_1/x_1')p(x_2/x_2') + p(x_1/x_2')p(x_2/ x_1')) \times% \nonumber\\
(p^i_1(x_1)p^i_2(x_2) + p^i_1(x_2)p^i_2(x_1))dx_1dx_2 \nonumber\\
= p^{i-}_1(x_1)p^{i-}_2(x_2) + p^{i-}_1(x_2)p^{i-}_2(x_1),
\end{align}
where $p^{i-}_1(x_1) = \int p(x_1/x_1')p^i_1(x_1')dx_1'$, i.e., the predicted pdf of prior pdf $p^i_1(.)$, and $p^{i-}_2(x_2)$ is the predicted pdf for the prior pdf $p^i_2(.)$. In H-FISST, the predicted MT-pdf for the $i^{th}$ component would simply be:
\begin{align}
p^-_i(\{x_1, x_2\}) = \int p(x_1/x_1')p(x_2/x_2') p^i_1(x_1')p^i_2(x_2')dx_1'dx_2', \label{FISST_prediction}
\end{align}
i.e., it is the same as the FISST pdf without the permutations of the arguments $\{x_1,x_2\}$. \\

Next, we turn to the update step. Suppose that we get the observation $\{z_1, z_2\}$. Again, we assume this purely for notational convenience and transparency of the treatment, it can easily be extended to more general observations at the expense of more notation. The FISST MT-likelihood function is then:
\begin{align}
p(\{z_1,z_2\}/\{x_1,x_2\}) %\nonumber\\
= p_D^2/ 2\{p(z_1/x_1)p(z_2/x_2) + p(z_1/x_2)p(z_2/x_1\}\nonumber\\
+ p_D(1-p_D)/2\{p(z_1/x_1)g(z_2) + p(z_2/x_1)g(z_1)\} \nonumber\\
+p_D(1-p_D)/2\{p(z_1/x_2)g(z_2) + p(z_2/x_2)g(z_1)\}%\nonumber\\
+(1-p_D)^2 \{g(z_1)g(z_2)\}, \label{FISST_likelihood}
\end{align}
where $p(z/x)$ is the single target likelihood, $g(z)$ is the probability of getting observation $z$ from clutter and $p_D$ is the probability of detection. The different terms in the likelihood function above correspond to the 7 different data associations possible given the two observations $\{z_1, z_2\}$ such as $(z_1\rightarrow T_1, z_2\rightarrow T_2)$, $(z_1 \rightarrow T_2, z_2 \rightarrow T_1)$, $(z_1 \rightarrow T_1, z_2 \rightarrow \mathcal{C})$, and so on. The factors of 1/2 in the first three terms and 1 in the fourth term of the likelihood equation above correspond to the ${m \choose k} k!$ normalization factor required in the MT-likelihood function: the 1/2 factor corresponding to $p_D^2$ term is ${2 \choose 2} 2!$, the 1/2 actor corresponding to the $p_D(1-p_D)$ term is ${2 \choose 1} 1!$ whereas the factor of 1 corresponding to $(1-p_D^2)$ is due to ${2\choose 0} 0!$. It may also be seen that due to the above normalization, $\int p(\{z_1, z_2\}/ \{x_1, x_2\}) dz_1dz_2 = 1$ for all $\{x_1, x_2\}$ as is generally true for a standard likelihood function in tracking/ filtering.
The updated FISST MT-pdf for the $i^{th}$ component is then given by the equation:
\begin{align}
p_i(\{x_1,x_2\}) = \frac{1}{\eta}\omega_i p(\{z_1,z_2\}/\{x_1,x_2\})p^-_i(\{x_1, x_2\}), \label{FISST_update}
\end{align}
where $\eta$ is a suitable normalization factor that will be evaluated below and is critical to understanding the structure of the FISST pdf. 
Consider the data association $(z_1\rightarrow T_1, z_2\rightarrow T_2)$, and call it the $i1^{th}$ association. Using Eqs. \ref{FISST_prediction}, \ref{FISST_likelihood} and \ref{FISST_update}, it may be seen that in the product in Eq. \ref{FISST_update}, this data association corresponds to the term:
\begin{align}
p_{i1}(\{x_1, x_2\})= \omega_ip_D^2/2 \{p(z_1/x_2)p(z_2/x_2)p^{i-}_1(x_1)p^{i-}_2(x_2) + %\nonumber\\
p(z_1/x_2)p(z_2/x_1) p^{i-}_1(x_2)p^{i-}_2(x_1)\}. \nonumber
\end{align}
It can be seen that the braced term in the expression above is nothing but the FISST MT-pdf:
\begin{align}
\eta_{i1} \{p^{i1}_1(x_1)p^{i1}_2(x_2) + p^{i1}_1(x_2) p^{i1}_2(x_1)\}, \mbox{where}\nonumber\\
\eta_{i1} = %\nonumber\\
\frac{1}{2!} \int \{p(z_1/x_1)p(z_2/x_2)p^{i-}_1(x_1)p^{i-}_2(x_2) + %\nonumber\\
p(z_1/x_2)p(z_2/x_1)p^{i-}_1(x_2)p^{i-}_2(x_1)\} dx_1dx_2\nonumber\\
= \int p(z_1/x_1) p(z_2/x_2) p^{i-}_1(x_1)p^{i-}_2(x_2)dx_1dx_2,
\end{align}
where
\begin{align}
p^{i1}_1(x_1) = \frac{p(z_1/ x_1)p^{i-}_1(x_1)}{\int p(z_1/x_1') p^{i-}_1(x_1')dx_1'},\nonumber\\
 p^{i1}_2(x_2) = \frac{p(z_2/ x_2)p^{i-}_2(x_2)}{\int p(z_2/x_2') p^{i-}_2(x_2')dx_2'},
 \end{align}
i.e., $p^{i1}_1(x_1)$ is simply the updated target 1 prior pdf  $p^{i-}_1(.)$ with observation $z_1$ and $p^{i1}_2(.)$ is the updated target 2 prior pdf $p^{i-}_2(.)$ with observation $z_2$. Similarly, the sub-components $p_{ij}(.)$ corresponding to the other possible data associations $ij$ may be found. Recall that $p_{ij}$ represent the transition probability from the $i^{th}$ parent to its $j^{th}$ child hypothesis.  Please note the distinction from the function $p_{ij}(\{X\})$ above which represents the MT-pdf underlying the $j^{th}$ hypothesis (we apologize for the notational ambiguity here but we wanted to be consistent with our HFISST derivation). Most importantly, after noting that $p_D^2/ 2 = p_{i1}$ in the H-FISST formulation, it may be seen that the normalization factor $\eta$ in Eq. \ref{FISST_update} can be written as:
\begin{align}
\eta = \sum_{i,j} \eta_{ij}p_{ij}\omega_i.
\end{align}
Thus, the ${i1}^{th}$ component of the FISST MT-pdf can be written as:
\begin{align}
\frac{\omega_i p_{i1} \eta_{i1}}{\sum_{i',j'}\eta_{i',j'} p_{i'j'} \omega_{i'}} \{p_1^{i1}(x_1)p_2^{i1}(x_2) + p_1^{i1}(x_2)p^{i1}_2(x_1)\}.
\end{align}
However, note that in the H-FISST framework $\eta_{i1} = l_{i1}$. Thus, in general, the ${ij}^{th}$ component of the FISST pdf may be written as:
\begin{align}
p_{ij}(\{x_1, x_2, \cdots x_n\}) = %\nonumber\\
\frac{\omega_i p_{ij} l_{ij}}{\sum_{i',j'} l_{i'j'}p_{i'j'}\omega_{i'}} \{\sum_{\bar{\sigma}} p^{ij}_1(x_{\sigma_1}) p^{ij}_2(x_{\sigma_2})\cdots p^{ij}_n(x_{\sigma_n})\},
\end{align}
where as before $\bar{\sigma} = \{\sigma_1, \sigma_2 \cdots \sigma_n\}$ represents all possible permutations of the numbers $\{1,2\cdots n\}$, and where $p^{ij}_1(x_1)p^{ij}_2(x_2)\cdots p^{ij}_n(x_n)$ is the updated MT-pdf that results from using the $j^{th}$ possible data association for component $i$. Note that this is precisely the update that is done in the H-FISST scheme modulo the set-theoretic representation in the FISST framework due to the interchangeability of the arguments $\{x_1, x_2\cdots x_n\}$. Hence, the above development shows that FISST and H-FISST recursions result in precisely the same MT-pdfs, and associated weights, modulo the representation of the underlying MT-pdfs in set theoretic terms in the FISST framework.\\
Further, note that as has been mentioned previously in the paper, the case of target birth and death results in there being more candidate children hypotheses ${ij}$, parameterized  through the transition probabilities $p_{ij}$. Thus, even in the case of target birth and death, the individual components/ hypotheses would have the same form as above in FISST, and thus, it follows that the H-FISST recursions and FISST recursions are the same.

\bibliographystyle{IEEEtran}
\bibliography{MAP_refs,FISST_refs,RFISST_Refs}

%%%%%%%%%%%%%%%%%%%%%%%%%%%%%%%%%%%%%%%%%%%%%%%%%%%%%%%%%%%%%%%%%%%%%%%%%%%%%%%%

\end{document}